\documentclass[a4paper,12pt]{article}

\usepackage{ucs}
\usepackage{amsfonts} 
\usepackage{amsmath} 
\usepackage{graphicx,subfigure}
\usepackage{caption} 
\usepackage{enumerate}
\usepackage{MnSymbol}
\usepackage{color}
\usepackage{xfrac}

\title{}
\author{}
\date{}

\newtheorem{theorem}{Theorem}[section]
\newtheorem{lemma}[theorem]{Lemma}

\newenvironment{definition}[1][Definition]{\begin{trivlist}
\item[\hskip \labelsep {\bfseries #1}]}{\end{trivlist}}

\newenvironment{@abssec}[1]{%
\if@twocolumn
\section*{#1}%
\else
\vspace{.05in}\footnotesize
\parindent .2in
{\upshape\bfseries #1. }\ignorespaces
\fi}
{\if@twocolumn\else\par\vspace{.1in}\fi}

\newcommand{\qed}{\nobreak \ifvmode \relax \else
\ifdim\lastskip<1.5em \hskip-\lastskip
\hskip1.5em plus0em minus0.5em \fi \nobreak
\vrule height0.75em width0.5em depth0.25em\fi}

\newcommand\keywordsname{Key words}
\newcommand\AMSname{AMS subject classifications}

\newcommand{\Q}{\mbox{WB}}
\newcommand{\WB}{\mbox{WB}}
\newcommand{\torus}{{\mathbb{T}^d}}

\newcommand{\R}{{\mathbb{R}}}

\newcommand{\Path }{}

\begin{document}
\title{Quantitative Quasiperiodicity}
\author{Suddhasattwa Das\footnotemark[1], \and Yoshitaka Saiki\footnotemark[2] \footnotemark[3] \footnotemark[5],
\and Evelyn Sander\footnotemark[4], \and James A Yorke\footnotemark[5]}
\footnotetext[1]{Department of Mathematics, University of Maryland, College Park}
\footnotetext[2]{Graduate School of Commerce and Management, Hitotsubashi University}
\footnotetext[3]{JST, PRESTO}
\footnotetext[4]{Department of Mathematical Sciences, George Mason University}
\footnotetext[5]{University of Maryland, College Park}
\date{\today}
\maketitle

\begin{abstract}{ 
The Birkhoff Ergodic Theorem concludes that time averages, i.e.,
Birkhoff averages, $\mbox{B}_N(f):=\Sigma_{n=0}^{N-1} f(x_n)/N$ of a function $f$ along
a length $N$ ergodic trajectory $(x_n)$  of a function $T$ converge to
the space average $\int f d\mu$, where $\mu$ is the unique invariant probability measure. 
Convergence of the time average to the space
average is slow. We use a modified average of $f(x_n)$ by giving
very small weights to the ``end'' terms when $n$ is near $0$ or $N-1$.
When $(x_n)$ is a trajectory on a quasiperiodic torus and $f$ and $T$
are $C^\infty$, our Weighted Birkhoff average (denoted $\Q_N(f)$) converges
``super'' fast to $\int f d\mu$ with respect to the number of iterates
$N$,  i.e. with error decaying faster than $N^{-m}$ for every
integer $m$. Our goal is to show that our Weighted Birkhoff average is a
powerful computational tool, and this paper illustrates its use for
several examples where the quasiperiodic set is one or two dimensional.
In particular, we compute rotation numbers and conjugacies (i.e. changes
of variables) and their Fourier series, often with 30-digit accuracy.
}\end{abstract}


\section{Introduction} \label{sect:intro}
Quasiperiodicity is a key type of observed dynamical behavior in a diverse set of applications. 
We say a map $T$ is ($d$-dimensionally) quasiperiodic (for $d\ge 1$)
if (i) $T: \torus \to \torus$ and (ii) each trajectory is dense in $\torus$ and (iii) there is a continuous choice of coordinates $\theta=(\theta_1,\cdots,\theta_d)$ and some $\rho=(\rho_1,\cdots,\rho_d)\in\torus$
for which the $T$ has the form
\begin{equation}\label{quasi}
T(\theta) = \theta+\rho \bmod1.
\end{equation}
Condition (ii) can be replaced by saying in dimension $d=1$ that $\rho$ is irrational and  in dimension $d>1$
that all of the coordinates $\rho_j$ of $\rho$ are irrational and they are ``independent'' over the reals; that is  if $a=(a_1,\cdots,a_d)$
is a vector of integers and $\sum_{j=1}^da_j\rho_j = 0$, then every $a_j =0$. We then say such a $\rho$ is {\bf irrational}.

Let $T$ be a $C^\infty$ quasiperiodic map. The quasiperiodicity persists for most small perturbations by the Kolmogorov-Arnold-Moser theory. 
We believe that quasiperiodicity is one of only three types of invariant sets with a dense trajectory that can occur in typical smooth maps.
The other two types are periodic sets and chaotic sets. 
See~\cite{sander:yorke:15} for the statement of our formal conjecture of this triumvirate. 
For example, quasiperiodicity occurs in a system of weakly coupled oscillators, in which there is an invariant smooth attracting torus in phase space with behavior that can be described exclusively by the phase angles of rotation of the system. 
Indeed, it is the property of the motion being described using only a set of phase angles that always characterizes quasiperiodic behavior. In a now classical set of papers, Newhouse, Ruelle, and Takens demonstrated a route to chaos through a region with quasiperiodic behavior, causing a surge in the study of the motion~\cite{newhouse:ruelle:takens:78}. 
There is active current interest in development of a systematic numerical and theoretical approach to bifurcation theory for quasiperiodic systems. 
Our goal in this paper is to present a numerical method for the fast calculation of the limit of Birkhoff averages in quasiperiodic systems, allowing us to compute various key quantities. 

If $f$ is integrable and the dynamical system is ergodic on the set in which the trajectory lives, then the Birkhoff Ergodic Theorem asserts that the Birkhoff average \boldmath $\mbox{B}_N$ \unboldmath defined as 
\begin{equation}\label{eq:B_N}
\mbox{B}_N(f):=\displaystyle\sum_{n=0}^{N-1} f(x_n)/N
\end{equation}
 of a function $f$ along an ergodic trajectory $(x_n)$ converges to the space average $\int f d\mu$ as $N\to\infty$ for $\mu$-almost every $x_0$, where $\mu$ is the unique invariant probability measure. 
In particular for quasiperiodic systems all trajectories with initial point in the ergodic set have the same limit of their Birkhoff averages. 
We develop a numerical technique for calculating the limit of such averages, where instead of weighting the terms $f(x_n)$ in the average equally, we weight the early and late terms of the set $\{0,\dots,N-1 \}$ much less than the terms with $n\sim N/2$ in the middle. 
That is, rather than using the equal weighting $(1/N)$ in the Birkhoff average, we use a weighting function $w(n/N)$.

\begin{figure}
\centering
\includegraphics[width = .5\textwidth] {\Path 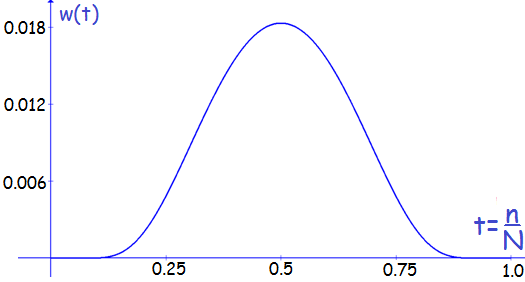}
\caption{\textbf{A \boldmath $C^\infty$\unboldmath~{\bf weighting function.}} The figure shows the graph of a the function $w(t)$ defined in Eq. \ref{eqn:weight}. This function plays the role of a temporal weighting for the Weighted Birkhoff average. That is, rather than using the equal weighting $(1/N)$ in the usual Birkhoff average, we use the non-uniform weight proportional to $w(n/N)$, where $n$ is the iterate number and $N$ is the orbit length. See Eqs. \ref{eqn:weight}, \ref{eqn:QN} for a description of this Weighted Birkhoff average.}

\label{fig:weight_exp}
\end{figure}

{\bf Weighted Birkhoff averaging method.}
A function $w:\R\to [0,\infty)$ will be called a 
\boldmath $C^\infty$\unboldmath~{\bf weighting function} 
if $w$ is infinitely differentiable and $w>0$ on $(0,1)$ and $= 0$ elsewhere. The example of such a function that we will use in this paper is in Eq. \ref{eqn:weight}  (Fig. \ref{fig:weight_exp}) defined as
\begin{equation}\label{eqn:weight}
w(t) :=\begin{cases}
\exp\left(\cfrac{-1}{t(1-t)}\right), & \mbox{for } t\in(0,1)\\
0, & \mbox{for } t\notin(0,1).
\end{cases}
\end{equation}

See Eq. \ref{eqn:wp} for a family of weighting functions $w^{[p]}(t)$ that converge even faster when many digits of precision are required.
This is an example of what is often referred to as ``window'' functions in spectral analysis or a ``bump'' function in the theory of partitions of unity.

For $d=1$, $\torus$ is a simple closed curve. For a continuous function $f$ and a $C^\infty$ quasiperiodic map $T$ on $\torus$, let $x_n\in\torus$ be such that $x_n = T(x_{n-1})$ for all $n\ge 1$. We define a {\bf Weighted Birkhoff ($\Q_N$) average} of $f$ as 
\begin{equation}\label{eqn:QN}
\Q_N(f)(x) :=\sum_{n=0}^{N-1} \hat{w}_{n,N}f(x_n),\mbox{ where }\hat{w}_{n,N}=\frac{w(n/N)}{\sum_{j=0}^{N-1}w(j/N)}.
\end{equation}

Note that $\Q_N(f)$ is indeed an average of the values $f(x_n)$ since $\sum_{n=0}^{N-1}\hat{w}_{n,N}=1$.

{\bf Main convergence result.} The main convergence result we are using is Theorem \ref{thm:A}. It is proved in \cite{Das-Yorke} and an outline of the proof is given here in Section~\ref{sect:nearly-rational}. 
We now state a special case of the theorem that avoids unnecessary terminology and states only the $C^\infty$ case. 

{\it
Assume $f$ and $T$ are $C^\infty$ 
and $w$ is a $C^\infty$ weighting function, 
then for almost every rotation number $\rho$ and for every positive integer $m$ there is a constant $C_m >0$ such that
}
\[\left| \Q_N(f)(\theta) - 
\int f(\theta) d\theta \right| \le C_m N^{-m}. 
\]

We refer to the above as {\bf super-fast} (super polynomial) {\bf convergence} or {\bf exponential convergence}. The above constant $C_m$ depends on (i) $w(t)$ and its first $m$ derivatives; (ii) the function $f(t)$ ; and (iii) the rotation number(s) of the quasiperiodic trajectory or more precisely, the small divisors arising out of the rotation vector. Our method of averaging does not give improved convergence results for chaotic systems. 

In \cite{Laskar99,Laskar93a,Laskar93b,Laskar03}, Laskar employs a Hanning data weighting function and the analogue of our $\sin^2$ weighting function  for his computations, which lead to the convergence of order $1/N^2$ or $1/N^4$, and $1/N^3$ respectively, (p. 136 in \cite{Laskar99})
where $N$ is the length of a orbit.  
Specifically, the weighting function $\cos^2(\pi x)$ for $x\in (-1,1)$ and averages over iterates from $-N$ to $N$.
He mentions the $C^{\infty}$ filter we use in Remark 2 of the appendix of \cite{Laskar99}, p. 146, though he does not use it.
There seems to be no advantage to using his lower order methods than  $C^{\infty}$ filter. In particular the programming of both is quite simple.
We will compare the two methods in Figs. \ref{fig:mid_circle2},  \ref {fig:2DMap_overview}(b), \ref {fig:circle_rotation}(b), and \ref {fig:2DMap_lyap}.

Other authors have considered related numerical methods (see Section \ref{Relatedmethods}), in particular~\cite{seara:villanueva:06, luque:Rot, luque:villanueva:14}, which we will compare to our approach when we introduce our averaging method in Section \ref{sec:methods}. 
See also~\cite{Simo1, Simo2, Durand:2002ug,Baesens:1990vj,Broer:1993uv,Vitolo:2011it,Hanssmann:2012jy, Sevryuk:2012dv,Broer:1990ip,Kuznetsov:2015dl}.
We announced some of the results presented here in \cite{DDSSSWY}.

{\bf The Babylonian problem of quasiperiodic rotation numbers.}
What constitutes a ``big-data'' problem depends on the speed of computation available. With this understanding, the first big-data problem was 2500+ years ago when the Babylonians computed the three periods of the moon from data on the position of the moon collected almost daily for many years. The moon's position through the fixed stars can be viewed as a quasiperiodic trajectory with $d=3$ and their problem was to compute the rotation numbers from such a trajectory, which they did with high accuracy. See
\cite{DSSYQR}.

{\bf Applications.}
We demonstrate our Weighted Birkhoff averaging method and its convergence rate by computing rotation numbers, conjugacies (i.e. changes
of variables), and their Fourier series in dimensions one and two. We will refer to a one-dimensional quasiperiodic curve as a \emph{curve}.

\begin{figure}
\centering
\subfigure[ ]{\includegraphics[width = .48\textwidth] {\Path 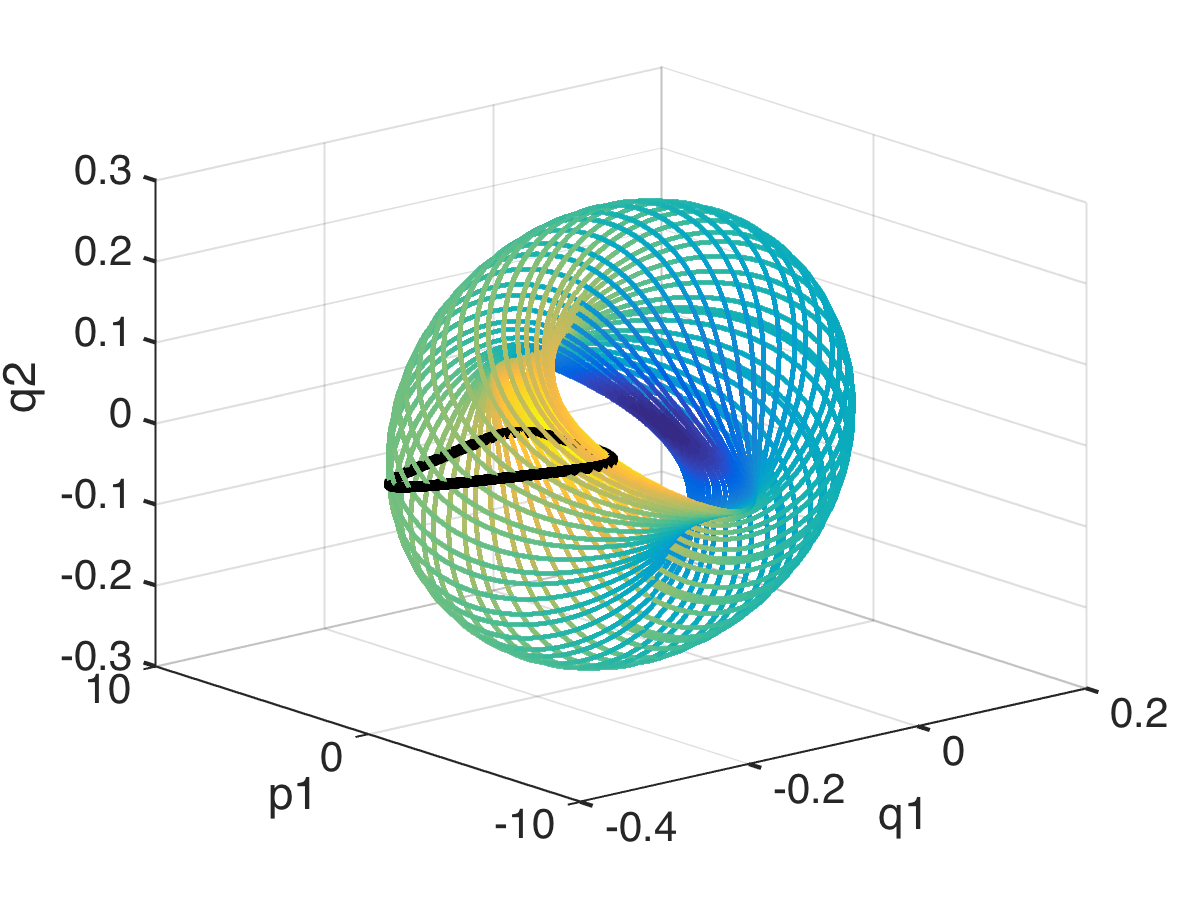}}
\subfigure[ ]{\includegraphics[width = .48\textwidth] {\Path 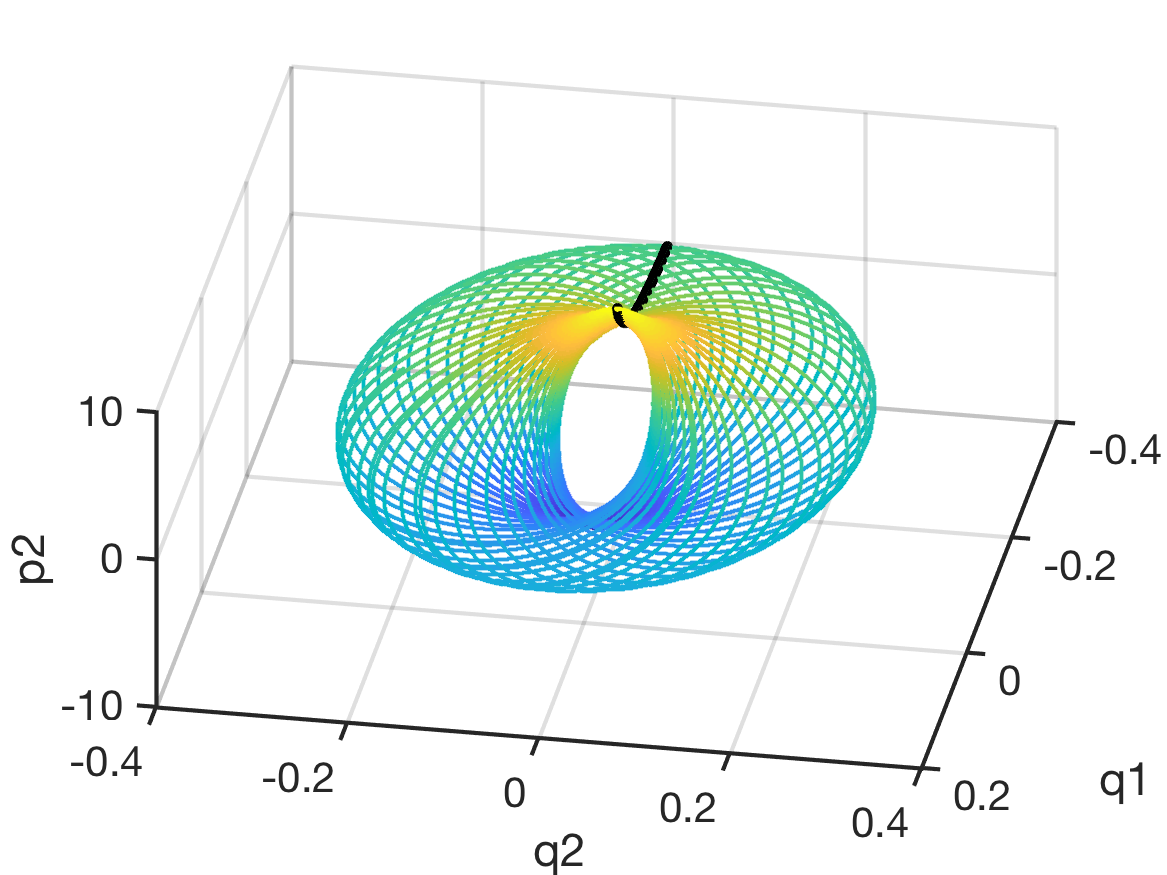}}

\subfigure[ ]{\includegraphics[width = .48\textwidth] {\Path 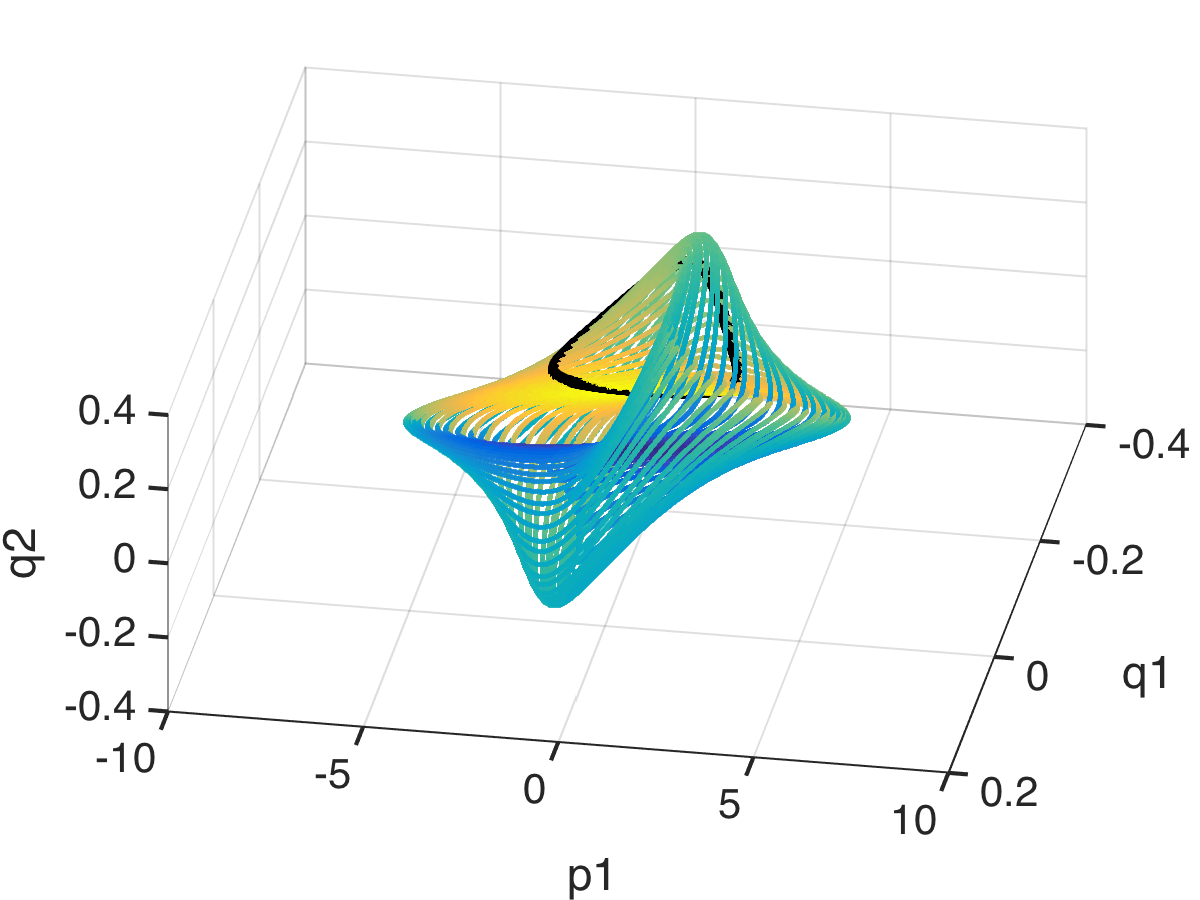}}
\subfigure[ ]{\includegraphics[width = .48\textwidth]{\Path 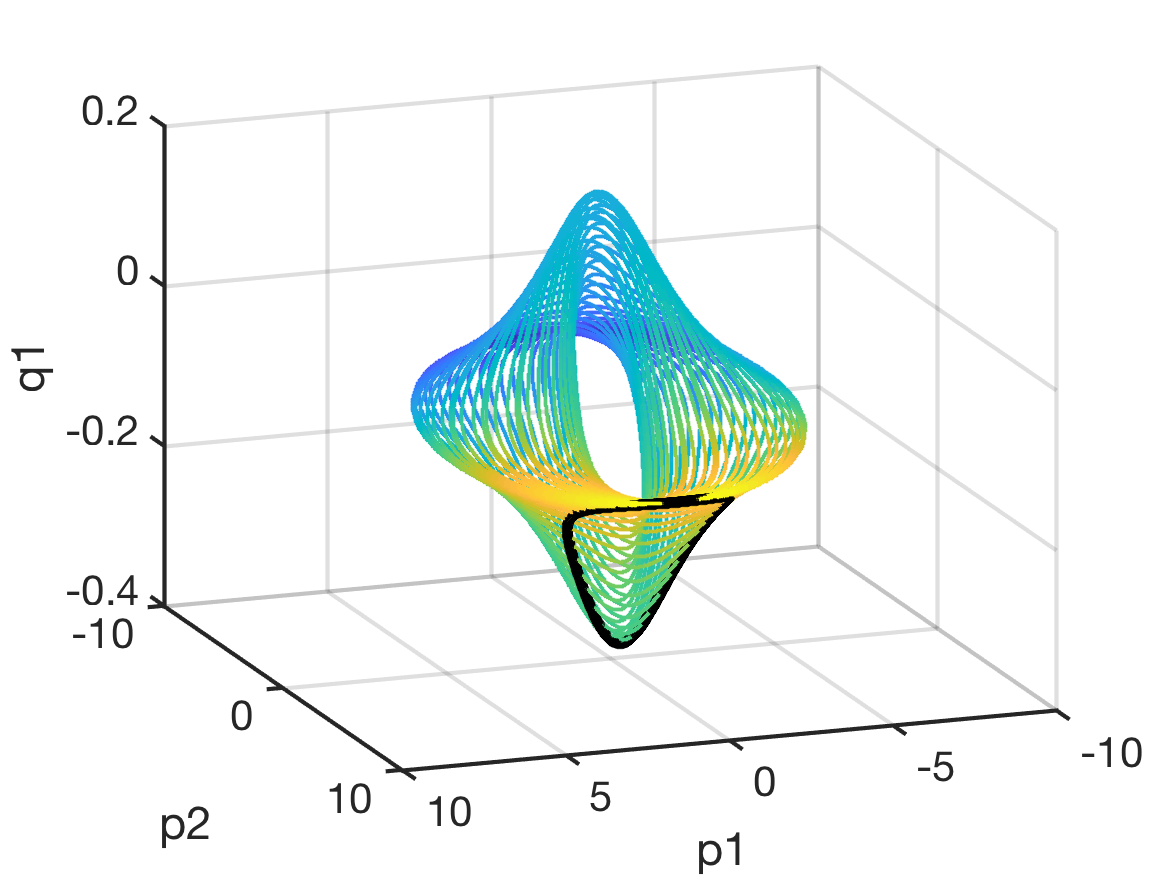}}
\caption{{\bf Torus flow for the R3BP.} 
This trajectory is a solution of Eq. \ref{eqn:ThreeBody}, shown as curve $B_1$ in Fig. \ref{fig:3B_global}.
All four views are of the same two-dimensional quasiperiodic torus lying in $\mathbb{R}^4$. Each picture consists of the same trajectory spiraling densely on this torus.  
We require four different views of this torus because the embedding into three dimensions gives a highly non-intuitive images. 
In all four panels, the color of the trajectory is the value of the variable $p_2$.
The black curve is the set of values of the Poincar{\'e} return map with $q_2=0$ for this flow torus.}
\label{fig:3B3}
\end{figure}

\begin{figure}
\centering
\subfigure[ ]{\includegraphics[width = .48\textwidth] {\Path 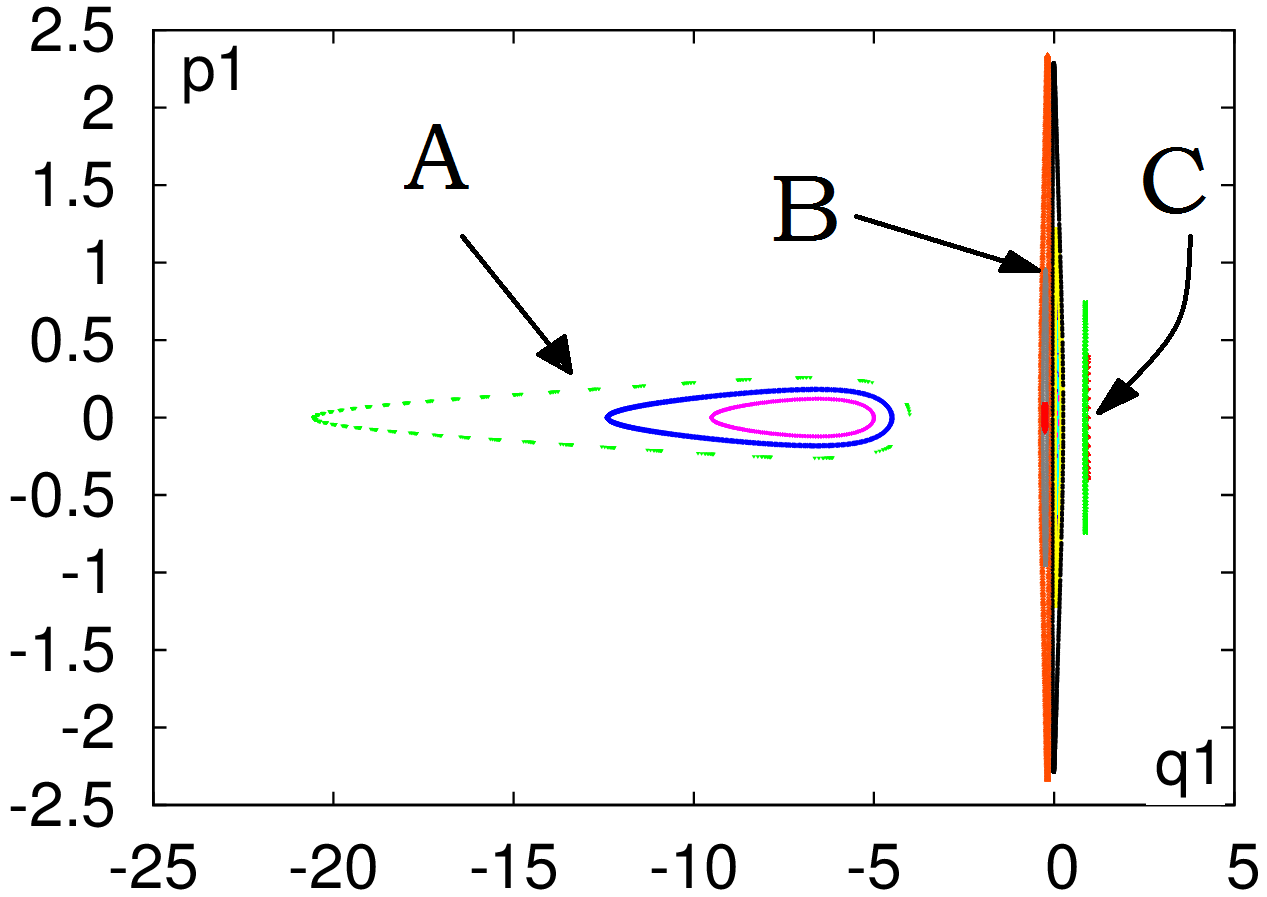}}
\subfigure[ ]{\includegraphics[width = .48\textwidth] {\Path 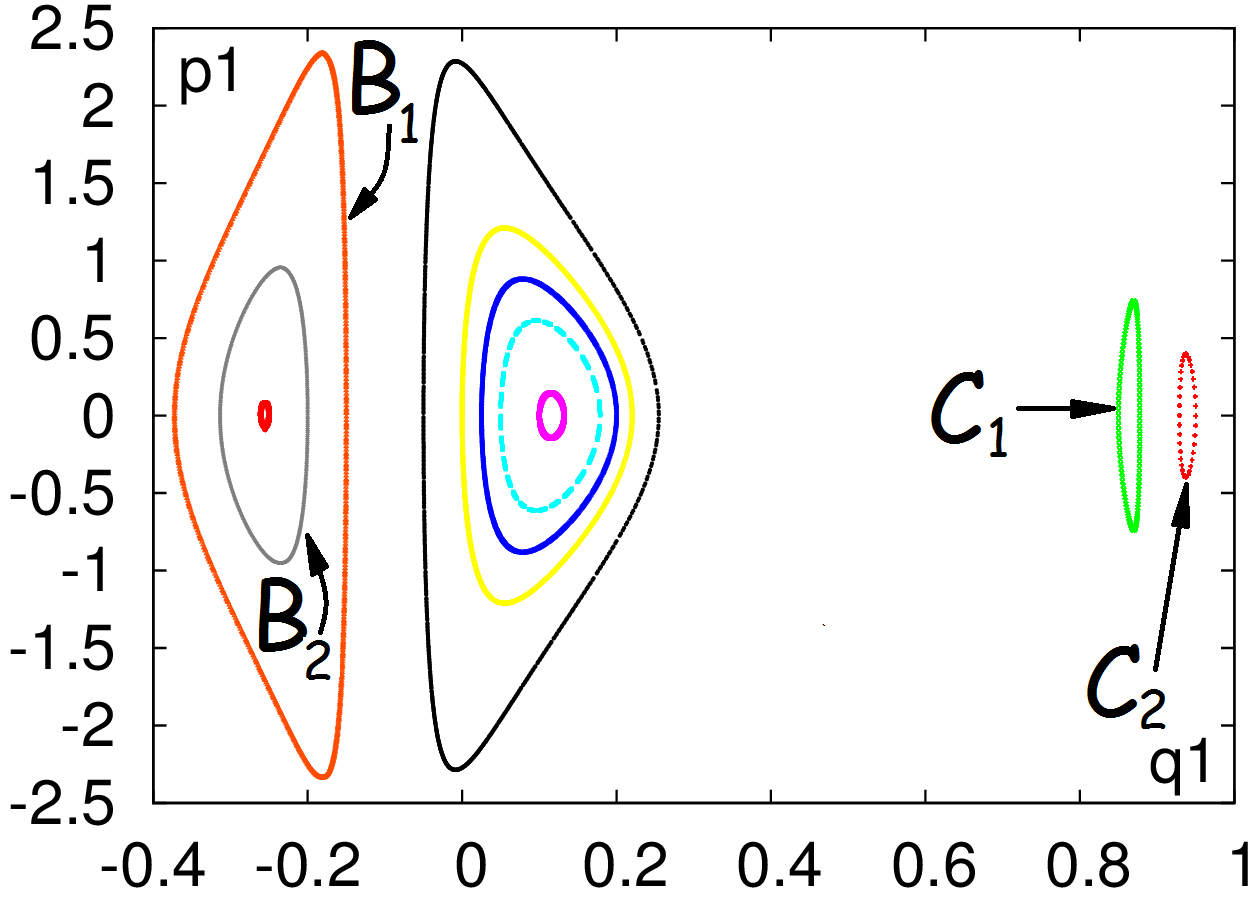}}
\subfigure[ ]{\includegraphics[width = .48\textwidth] {\Path 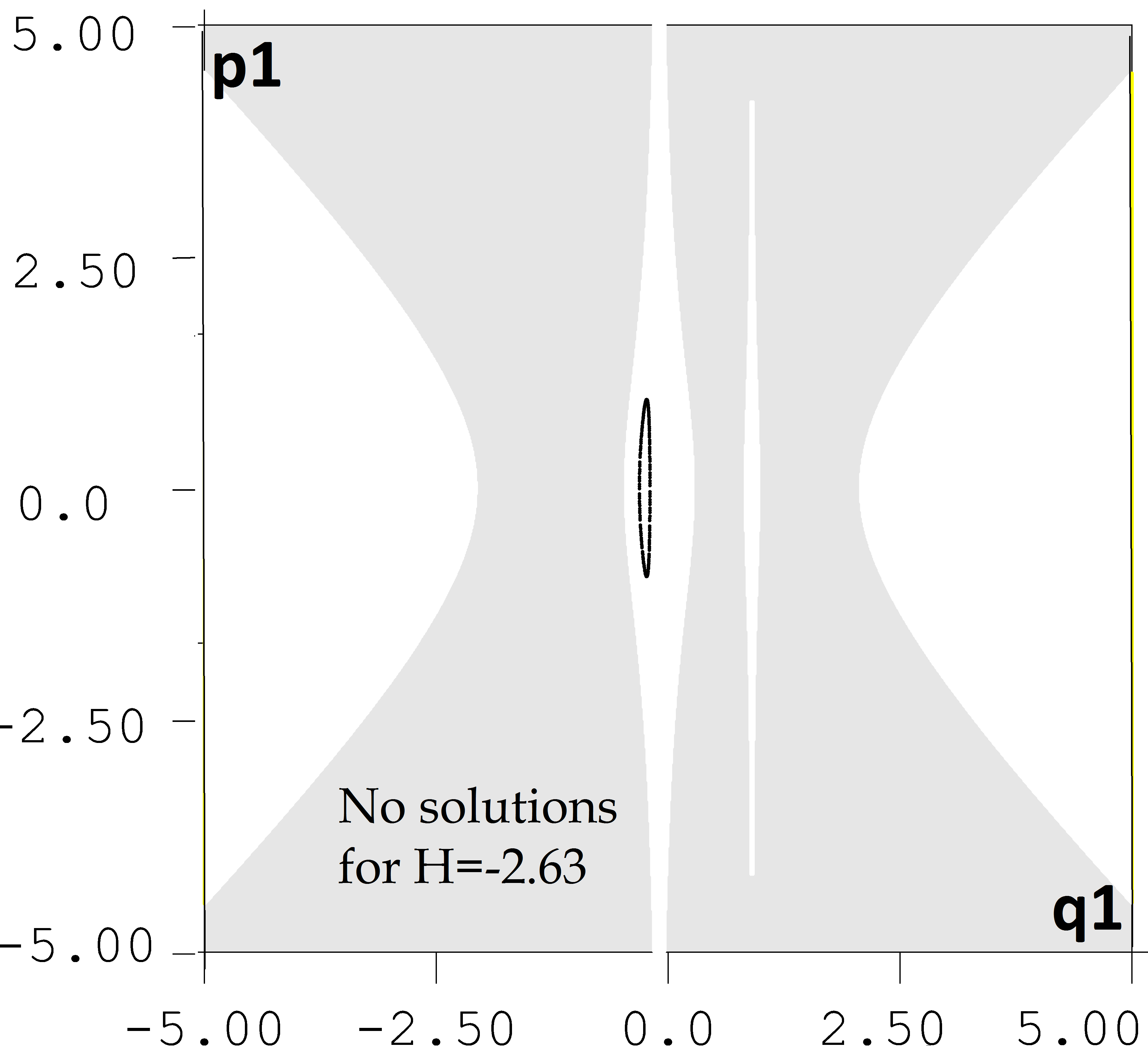}}
\caption{ {\bf Poincar{\'e}-return map for R3BP.} \label{fig:3B1}
Panels (a) and (b) show a projection to the $q_1-p_1$ plane of solutions to the R3BP  in Eq.~\ref{eqn:ThreeBody}. The value of the Hamiltonian $H$ for all the curves shown is the same and $H\approx -2.63$. They show various quasiperiodic trajectories on the Poincar{\'e} section corresponding to $q_2=0$. Each time the flow hits $q_2=0$   and $dq_2/dt>0$, we plot $(q_1,p_1)$. 
Note that the planet is fixed at the point $(-0.1,0)$ and the moon at $(0.9,0)$. Thus some trajectories orbit both the planet-moon system and some orbit only the planet or only the moon. Each trajectory shown is a quasiperiodic curve. 
In Panel (c), white indicates the region of the $(q_1,p_1)$ plane where the Poincar{\'e} return map is defined for $H\approx -2.63$, and gray indicates the region of the plane where the Poincar{\'e} return map is not defined. Panel (c) also shows the trajectory which corresponds to the curve $B_1$. 
}
\label{fig:3B_global}
\end{figure}
We start by describing our results for a key example of quasiperiodicity: the (circular, planar) restricted three-body problem ({\bf R3BP}). This is an idealized model of the motion of a planet, a large moon, and an asteroid governed by Newtonian mechanics, in a model studied by Poincar{\'e}~\cite{ThreeBody2,ThreeBody1}. In particular, we consider a planar three-body problem consisting of two massive bodies (``planet" and ``moon") moving in circles about their center of mass and a third body (``asteroid") whose mass is infinitesimal, having no effect on the dynamics of the other two.

We assume that the moon has mass $\mu$ and the planet mass is $1-\mu$ where $\mu = 0.1$, and
writing equations in rotating coordinates around the center of mass. Thus the planet remains fixed at $(q_1,p_1) = (-0.1,0)$, and the moon is fixed at $(q_2,p_2)=(0.9,0)$. In these coordinates, the satellite's location and velocity are given by the {\em generalized position vector} $(q_1,q_2)$ and {\em generalized velocity vector} $(p_1,p_2)$. The equations of motion are as follows.  
%
%
\begin{equation}
{\displaystyle
\begin{array}{rcl}\label{eqn:ThreeBody}
\cfrac{dq_1}{dt} &=& p_1+q_2, \\
\cfrac{dq_2}{dt} &=& p_2-q_1, \\
\cfrac{dp_1}{dt} &=& p_2- \mu\cfrac{q_1-1+\mu} {d_{moon}^{\ 3}} -(1-\mu)\cfrac{q_1+ \mu}{ d_{planet}^{\ 3}},\\
\cfrac{dp_2}{dt} &=& -p_1
-\mu\cfrac{ q_2} {d_{moon}^{\ 3}}
-(1-\mu)\cfrac{q_2}{ d_{planet}^{\ 3}},
\end{array}
}
\end{equation}

where 
\[ d_{moon}^{\ 2} = (q_1-1+\mu)^2+q_2^2\] \[d_{planet}^{\ 2} = (q_1+\mu)^2+q_2^2.\]

The following function $H$ is a Hamiltonian for this system
\begin{equation}\label{eqn:Hamiltonian}
H=\frac{p_1^2+p_2^2}{2}
+\left(p_1q_2-p_2q_1\right)
+\left(-\frac{1-\mu}{d_{planet}} - \frac{\mu}{d_{moon}}\right), 
\end{equation}
where $p_1=\dot{q_1}-q_2$ and $p_2=\dot{q_2}+q_1$ (see \cite{ThreeBody3} p.59 Eqs. 63-66).

The three terms in Eq. \ref{eqn:Hamiltonian} are resp. the kinetic energy, angular moment, and the potential.
For fixed $H$, Poincar{\'e} reduced this problem to the study of the Poincar{\'e} return map for a fixed value of $H$, only considering 
a discrete trajectory of the values of $(q_1,p_1)$ on the section $q_2=0$ and $dq_2/dt>0$. 
Thus we consider a map in two dimensions rather than a flow in four dimensions. 
Fig. \ref{fig:3B3} shows one possible motion of the asteroid for the full flow. The 
orbit is spiraling on a torus. The black curve shows the corresponding trajectory on the 
Poincar{\'e} return map. Fig. \ref{fig:3B1} shows the Poincar{\'e} return map 
for the asteroid for a variety of starting points. 
A variety of orbits are shown, most of which are quasiperiodic invariant curves. 
An exception is trajectory A in Fig. \ref{fig:3B_global}(a), which is an invariant recurrent set consisting of $42$
curves. Each curve is an invariant quasiperiodic curve under the 42$^{nd}$ iterate of the map. 

Our paper proceeds as follows: In Section~\ref{sec:Quasi}, we give the formal
definition of quasiperiodicity, rotation number, and the conjugacy map to the rigid rotation. 
In Section~\ref{sec:methods}, we describe our
numerical technique in detail. We illustrate our Weighted Birkhoff averaging method for 
a series of four examples, including an example of a two-dimensionally quasiperiodic map. 
In all cases, we get fast convergence and are in most cases able to 
give results with about 30-digit accuracy. Section \ref{sect:nearly-rational} describes what happens when a rotation number is unusually well approximated by a fraction with small denominator. In such cases we elliptically say the rotation number is ``nearly rational''.
Finally, Section~\ref{sec:conclude} contains our concluding remarks. 
\section{Quasiperiodicity}\label{sec:Quasi}

In the introduction, we described quasiperiodic motion as motion that
could be fully understood through a set of angles of rotation. We now
formalize that idea in the following definition.

\textbf{Quasiperiodicity.} For a dimension $d \ge 1$, let $\rho = (\rho_1,\rho_2,\dots, \rho_d)$ be a 
vector whose coordinates are irrational and are independent over the integers
(see Eq. \ref{quasi}). The following map $T_{\rho}: \torus \to \torus$ is called a {\bf rigid rotation}:
\begin{equation}\label{eqn:Rigid_rotation}
T_{\rho}(\theta) = \theta +  \rho \bmod 1, \mbox{ where ``mod'' is applied to each coordinate.} 
\end{equation}
A rigid rotation is the simplest, albeit least interesting example of a map with quasiperiodicity. Since $T_{\rho}$ gives the same values on opposite sides
of the unit cube, we identify the sides and refer to the domain of the rigid rotation as a {\bf curve} in one dimension and a {\bf $d$-torus} in  dimension $d>1$. In this paper, we will sometimes refer to the curve as a 1-torus. 
We define a map $T$ on an ambient space $\R^D$ to be {\bf quasiperiodic} if either $T$ or some iterate $T^k$ is topologically conjugate to a rigid rotation. 
We will assume $k=1$ in the rest of this description. That is, a map $T$ is quasiperiodic if there is a rigid rotation map $T_{\rho}$ and an invertible {\bf conjugacy map} (i.e., change of coordinates) 
$h: \torus \to \R^D$ such that
\begin{equation}
T(h(\theta)) = h(T_{\rho}(\theta)). \label{eqn:conj}
\end{equation} 
A flow has quasiperiodic behavior if its associated Poincar{\'e} return map has quasiperiodic behavior.

For an invertible map $T$ to be quasiperiodic on a curve $C$, it is necessary and sufficient that  $T$ has a dense trajectory, as shown in \cite{TransCircleHomeo}. In general, a one-time differentiable invertible map on a curve without periodic points may not be quasiperiodic. However, if we assume that the map $T$ and the curve $C$ are twice continuously differentiable, then Denjoy~\cite{VanKampen} showed that these conditions are both necessary and sufficient. Furthermore, clearly any
rigid irrational rotation map is a real analytic map, but even if we
assume that a quasiperiodic function is analytic, Arnold showed that the conjugacy map $h$
may only be continuous for some atypical rotation number. However, Herman (see \cite{HermanSeminal}) proved that for 
homeomorphisms on a circle, most conjugacy maps $h$ are analytic. 
Yamaguchi and Tanikawa~\cite{StdMap1} and Hunt, Khanin, Sinai and Yorke~\cite{HKSY} show that the critical KAM curve may not be $C^2$.

{\bf Diophantine rotations.} An irrational vector $  \rho \in\mathbb{R}^d$ is said to be {\bf Diophantine} if for some $\beta>0$ it is {\bf Diophantine of class} $\beta$ (see \cite{HermanSeminal}, Definition 3.1), which means there exists $C_\beta>0$ such that for every $ k \in\mathbb{Z}^d$, $ k \neq 0$ and every $n\in\mathbb{Z}$,
\begin{equation}\label{eqn:Dioph}
| k \cdot \rho -n|\geq\frac{C_\beta}{\|k\|^{d+\beta}.}
\end{equation}
We conjecture that {\bf if the map is analytic, then almost every quasiperiodic torus having a rotation number that is Diophantine (i.e., far from rational) is real analytic.}

\begin{figure}
\centering
\subfigure[ ]{\includegraphics[width = .4\textwidth] {\Path 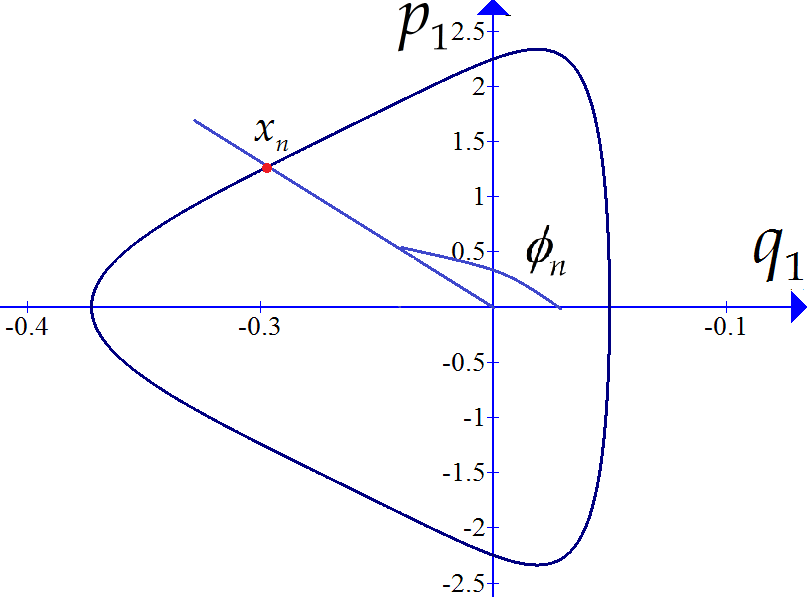}}
\subfigure[ ]{\includegraphics[height=.2\textwidth] {\Path 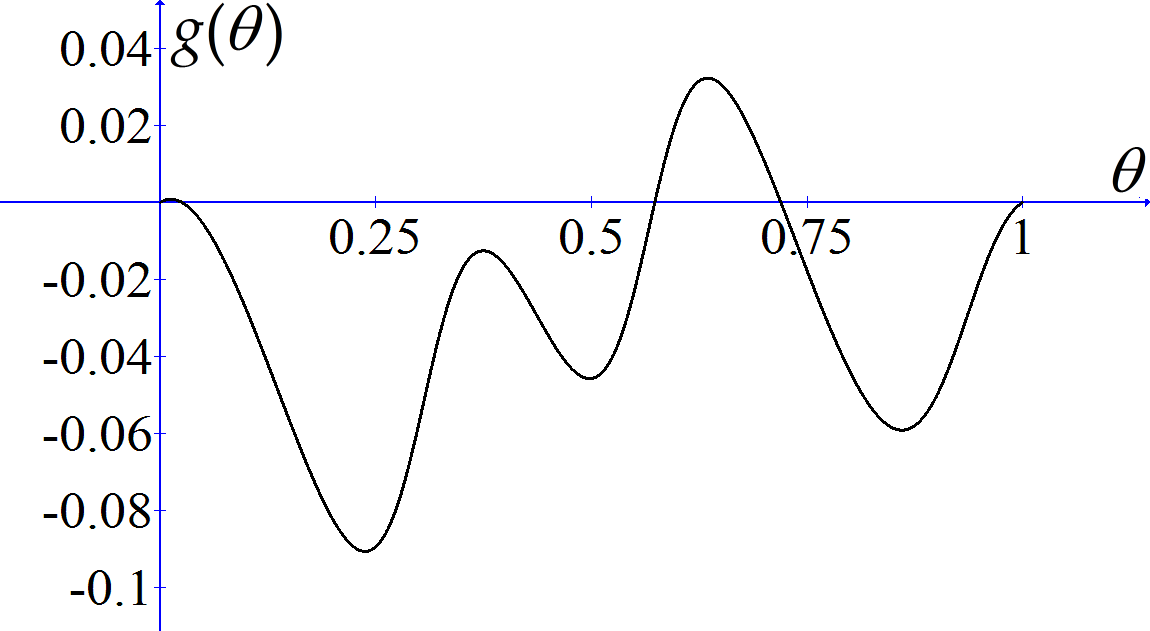}}
\subfigure[ ]{\includegraphics[height= .4\textwidth] {\Path 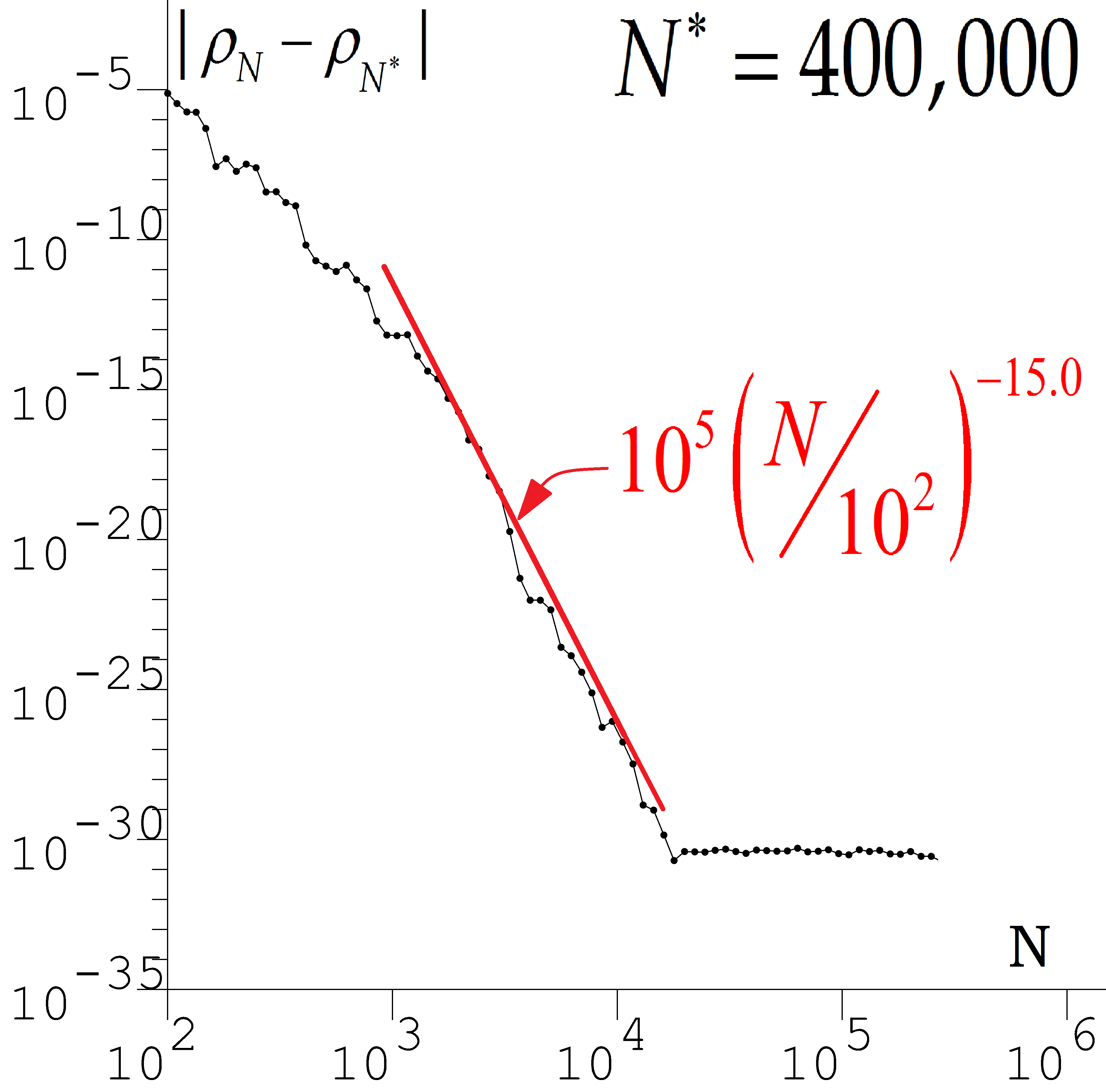}}
\subfigure[ ]{\includegraphics[height= .4\textwidth] {\Path 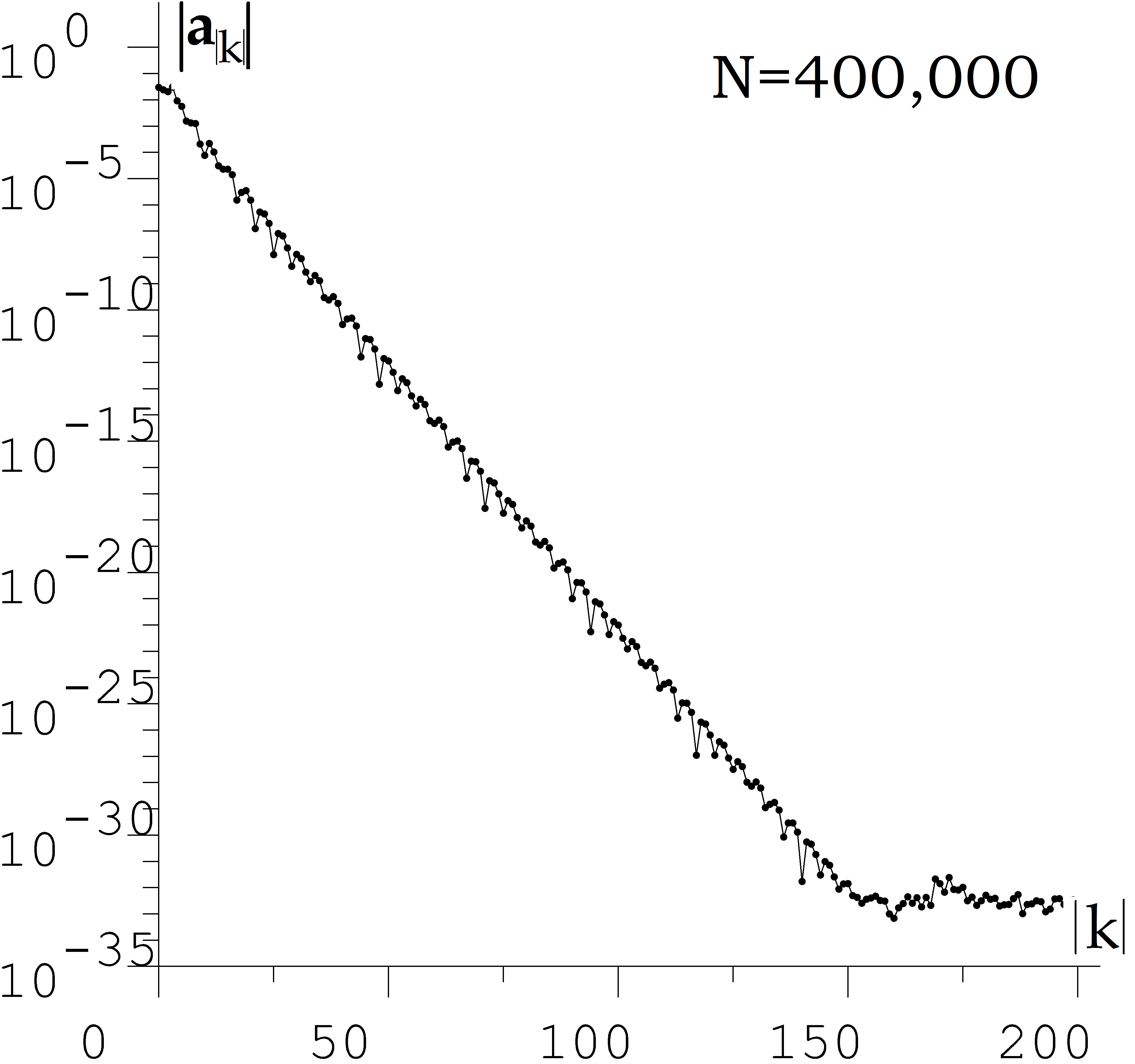}}
\caption{\textbf{Quasiperiodicity for the R3BP.} For the quasiperiodic curve $B_1$ from Fig. \ref{fig:3B1}, panel (a) shows the curve $B_1$, along with the projection mapping from $x_n$, the $n^{th}$ iterate in the $(q_1,p_1)$ coordinates,  to 
angular coordinate $\phi_n$. It is parameterized by coordinates $\phi \in S^1 \equiv [0,1)$. Panel (b) depicts the periodic part $g(\theta)$ of the one-dimensional conjugacy map (Eq. \ref{eqn:g_conjugacy}) between the quasiperiodic behavior and rigid rotation by $\rho$. 
The curve represents the true Fourier series up to 30 decimal digits. 
Panel (c) shows the convergence rate of the error in the rotation number $\rho_N$ as a function of the number of iterates $N$. 
The ``error'' is the difference $|\rho_N-\rho_{N^*}|$, where $N^*=400,000$ is large enough so that $\rho_N$ appears to have converged. 
The straight line is an upper bound for the curve and its exponent indicates the rate of convergence.
Panel (d) shows the norm of the Fourier coefficients of the conjugacy as a function of index. 
This exponential decay indicates that the conjugacy function is analytic, up to numerical precision.  
The step size used for the $8^{th}$ order Runge-Kutta scheme is $2 \times 10^{-5}$.}
\label{fig:results_3B}
\end{figure}

\textbf{Assigning angular coordinates.} Let $(x_n)$ be the forward orbit under $T$, and $( \theta_n )$ the forward orbit under $T_{\rho}$. That is, $x_{n+1} = T(x_n)$ and $\theta_{n+1} = \theta_n + \rho \pmod 1$.  

In all of the $d=1$ examples discussed here, the quasiperiodic curves are simple, closed, convex curves, and hence, angular coordinates can be obtained using polar coordinates $(\phi,r)$  with respect to a fixed and pre-defined center, where $r$ is uniquely determined by $\phi$. 
More generally we may have an image of a quasiperiodic curve or torus that is not an embedding as a closed, convex curve. Such is not needed, since we have a method based on the Takens delay coordinate maps that solves the problem.
This general method of obtaining rotation numbers from the image of a quasiperiodic curve is announced in~\cite{DDSSSWY} with a complete description in
\cite{DSSYQR}.
 
Because of our nice embeddings the orbits $(\theta_n=n\rho)$, $(x_n)$ and $(\phi_n)$ are all conjugate and there is a continuous map $V:S^1\to S^1$ such that
\begin{equation}\label{eqn:conjugacy1D}
x_n = h(\theta_n),\quad \phi_n = V(\theta_n),\quad \theta_n=n\rho\bmod 1~\mbox{ for all } n=0,1,2,3,\ldots 
\end{equation}
Since $V$ is invertible, the following map is periodic.
\begin{equation}\label{eqn:g_conjugacy}
g(\theta) := V(\theta)-\theta
\end{equation}
For example, in Figure~\ref{fig:results_3B}(a), the angular coordinate $\phi$ for the quasiperiodic curve $B_1$ of the restricted three-body problem  is measured from a point $(q_1,p_1) =(-0.2,0)$.
\section{Weighted Birkhoff averaging $\Q_N$ and its applications}\label{sec:methods}
As mentioned in Section~\ref{sect:intro}, our approach is to modify the regular Birkhoff average in Eq. \ref{eq:B_N} using the weighting function $w$ in Eq. \ref{eqn:weight} (cf. Fig. \ref{fig:weight_exp})  so that for quasiperiodic dynamical systems, the Weighted Birkhoff average in Eq. \ref{eqn:QN} convergences much faster to the same limit $\int f(\theta)d\theta$. We will formalize this claim using the main result from the companion paper \cite{Das-Yorke}.
\begin{definition}
A function $w:\mathbb{R}\to [0,\infty)$ is said to be a \boldmath $C^\infty$ \unboldmath \textbf{bump function} if $w$ is $C^\infty$ and
 the support of $w$ is $[0,1]$ and $\int_\mathbb{R}w(x)dx \ne 0$ and $w$ and all of its derivatives vanish at $0$ and $1$. 
\end{definition} 
\begin{theorem}[Theorems 1.1, 3.1 in \cite{Das-Yorke}]\label{thm:A}
For $r \in[1,\infty]$, let $X$ be a $C^r$ manifold and $T:X \to X$ be a $C^r$ map which is $d$-dimensionally quasiperiodic on an invariant set $X_0\subseteq X$, with invariant probability measure $\mu$ and a rotation vector of Diophantine class $\beta (>0)$. Let $f:X \to \mathbb{R}^k$ be a $C^{r}$ map. 
Let  $m>1$ be an integer such that $r\geq d+m(d+\beta)$, and let $w:\R \to \R$ be a $C^{m}$ bump function. Then there is a constant $C_m$ depending upon $w,f,m,M,$ and $\beta$ but independent of $x_0\in X_0$ such that
\begin{equation}\label{eqn:Mm2}
\left|(\WB_Nf)(x_0)- \int_{X_0}\!f~d\mu\right|  \le C_m N^{-m} ,
\end{equation}
In particular, if $r=\infty$ and $w$ is a $C^\infty$ bump function, then Eq. \ref{eqn:Mm2} holds for every $m\in\mathbb{N}$.
\end{theorem}
Theorem \ref{thm:A} is proved in Section \ref{sect:nearly-rational}
in a way that lets us determine what happens when an irrational rotation number is near a rational number.

{\bf Diophantine rotation numbers.} The assumption on the rotation numbers being Diophantine means that the rotation numbers cannot be closely approximated by rational numbers with small denominators. For numbers which are not Diophantine, the trajectories have ergodic properties ``close'' to periodic orbits and therefore, do not converge slower than the rate $N^{-k}$ for some $k$.

{\bf Robustness of assumptions.} It is well known that for every $\beta>0$ the set of Diophantine vectors of class $\beta$ have full Lebesgue measure in $\mathbb{R}^d$ (see for example, \cite{HermanSeminal}, 4.1). Thus, the assumption of the rotation number being Diophantine is robust in a measure theoretic sense, i.e., in physical experiments, the rotation number will be  Diophantine with probability 1.


\subsection{Computing a rotation number or rotation vector}\label{sec:Method_Rot}
We now show how to apply this averaging method in computation. We observe that $N$ must generally be larger for $\mathbb{T}^2$ than for $\mathbb{T}^1$ to get a 30-digit accuracy.

According to the definition of ``quasiperiodicity'', a quasiperiodic map is conjugate to a rigid rotation of the form Eq. \ref{eqn:Rigid_rotation} with rotation vector $\rho$. However, this vector $\rho$ is not unique. When the dimension $d=1$, there are two choices of $\rho$ depending on whether you move clockwise or counterclockwise around the circle, and for $d>1$, it is shown in \cite{DSSYQR} that there is a set of choices for $\rho$ which are related to each other by unimodular transformations of $\mathbb{R}^d$ and are dense in $\torus$. Our goal will be to find any one of these equivalent rotation vectors.
The rotation vector of a quasiperiodic trajectory $(y_n)_{n=0}^{N-1}$ can be viewed as the average rotation traversed by the sequence of $d$-vectors
$(T(\theta_n) - \theta_n) $, which by Eq. \ref{eqn:conjugacy1D} is equivalent to the average of the angular increments  $(\phi_{n+1}-\phi_n)$.  To make the notion of the angular distance between $\phi_{n}$ and $\phi_{n+1}$ consistent, we must make this choice continuously. By conjugating with $(hV^{-1})$ from Eq. \ref{eqn:conjugacy1D} if necessary, we may assume that the map $T$ is in angular coordinates on $\torus$. 

We then associate with $T$ a continuous map $\tilde{T}$ on the full Euclidean space $\mathbb{R}^d$ such that 
\begin{displaymath}\label{eqn:q_as_cover}
	\tilde{T}(z) \pmod 1 = T(z \pmod 1).
\end{displaymath}
The map $	\tilde{T}$ is called a {\bf lift} of $T$, and $z\in \mathbb R^d$ is a lift of $\phi\in\torus$. 

Since $T$ is invertible, the map  $\tilde{T}(z)-z$ has period one in every coordinate direction. For example, the rigid rotation $T(\phi)  = \phi + \rho \pmod 1,$ for $\theta$ in $\torus$ has a corresponding lift map $	\tilde{T}(z) = z + \rho$. 
Of course if $\rho$ was $\sqrt2$, we would have the same map $T$ as for $\rho = \sqrt 2 - 1$ so we define $\tilde T$ using rotation numbers are in $(0,1)$.
Using the lift, we now give a formula for the rotation vector for the trajectory $(y_n)$ starting at $y_0$: 
\begin{equation}\label{eqn:rot_num_defn}
\rho(T):=\lim_{N\rightarrow\infty} \frac{1}{N} \sum_{n=0}^{N-1} \left( 	\tilde{T}(z_n)-z_n \right) .
\end{equation}

This average converges slowly as $N\to\infty$, with order of at most $1/N$. However, since Eq. \ref{eqn:rot_num_defn} can be written as a Birkhoff average by writing $f(z_n)=\tilde{T}(z_n)-z_{n}$, we can apply our method to this function. That is, let $(z_n)_{n=0}^{N-1}$ be an orbit for $	\tilde{T}$.
Our approximation of $\rho$ is given by the Weighted Birkhoff average of $f$,
\begin{equation}\label{eqn:weighted_rho}
\WB_N (z_{n+1}-z_n) := \sum\limits_{n=0}^{N-1} \hat{w}_{n,N} (z_{n+1}-z_n)\rightarrow \rho \mbox{ as } N \to\infty.
\end{equation}

\subsection{Convergence rate of the Weighted Birkhoff average  $\Q_N$ } 
\begin{figure}
\centering
\includegraphics[height=.4\textwidth]{\Path 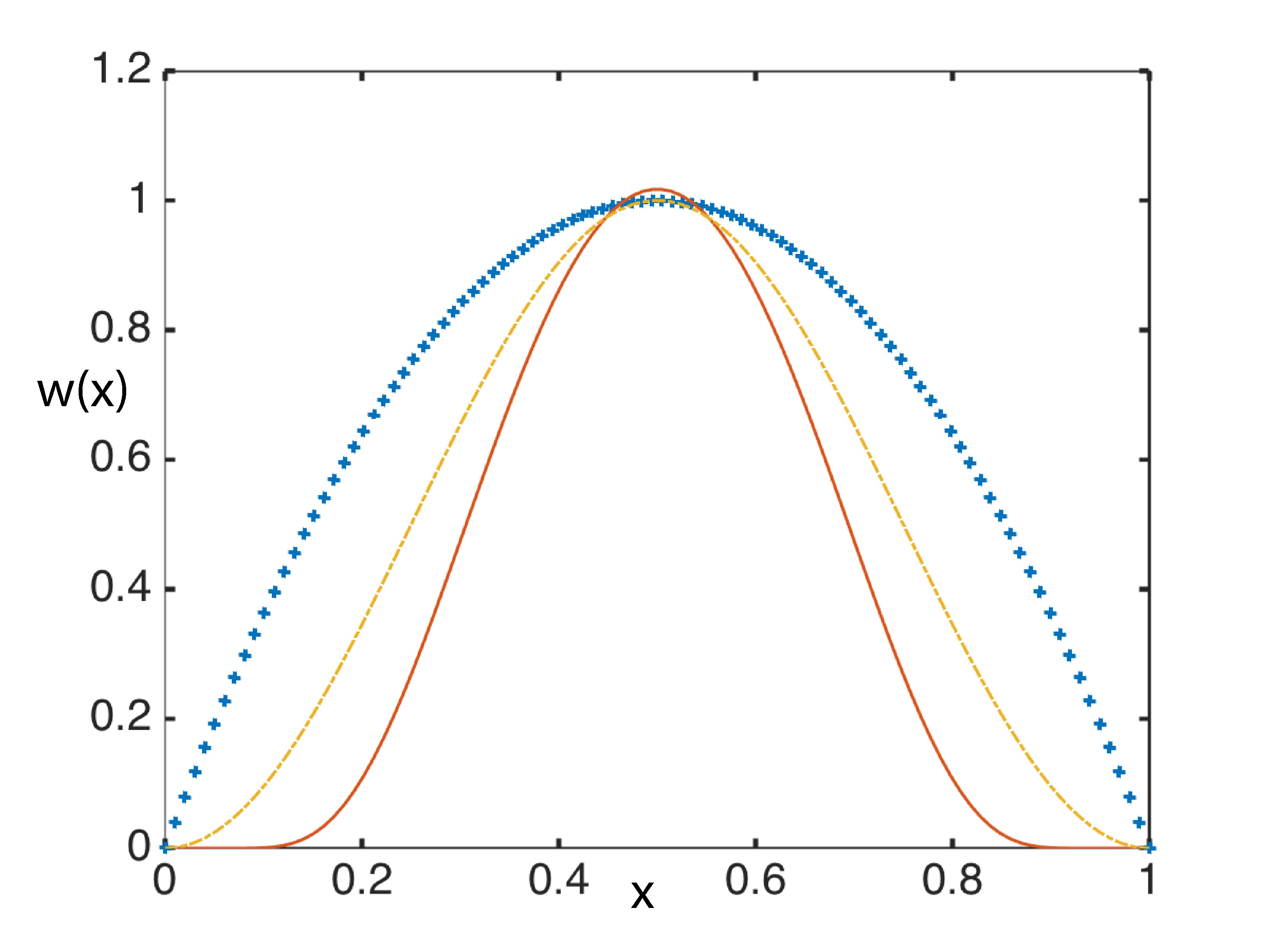}
\caption{\textbf{Variations of different weighting functions}. 
This is a plot of three non-constant weighting functions from Eq. \ref{eqn:weights} and below $w_{quad}$ (top), $w_{\sin^2}$ (second),
$w_{\exp}$ (lowest). Since only the shape matters, they have been 
rescaled so that each has a peak of approximately $1.0$. 
}
\label{fig:mid_circle}
\end{figure}

In order to illustrate the speed of convergence of our Weighted Birkhoff average $\Q_N$ as $N\rightarrow\infty$, we introduce four different possible choices for the weighting function $w$, depicted in Fig.~\ref{fig:mid_circle}, and compare the convergence results for computing the rotation number for each of these choices of $w$. 
\begin{eqnarray}
w_{equal}(t) & := & 1\mbox{ (Birkhoff's choice)}  \label{eqn:weights} \\
w_{quad}(t) & := & t (1-t)\nonumber  \\
w_{\sin^2}(t) & := & \sin^2(\pi t) \nonumber \\
w^{[1]}(t) := w_{\exp}(t) & := 
& \exp\left(\frac{-1}{t(1-t)}   \right). \nonumber 
\end{eqnarray}
The weighting functions described above are defined to be $0$ outside $(0,1)$ and equal to the specified function inside $(0,1)$. 
Recall that the last function in the list, $w_{\exp}$, is the function used in our calculations. 
It is the only one in the list that is $C^\infty$. The others fail to be $C^\infty$ at $t=0$ and $t=1$.
A family of $C^\infty$ weighting functions can be defined for $p\ge 1$ as 
\begin{equation}\label{eqn:wp}
w^{[p]}(t) := \exp\left(\frac
{-1}{t^p(1-t)^{p}}\right)
\mbox{ for } t\in(0,1)
\end{equation}
  and $=0$ elsewhere.  
This paper mainly uses $w= w^{[1]}$. The function $w^{[2]}$ results in an averaging method which converges noticeably faster than $w^{[1]}$ when using 30-digit precision, but not in 15-digit precision. We will write the Weighted Birkhoff averages as
$ \Q_N^{[1]}$ (or just $\Q_N$) and $\Q_N^{[2]}$ when using
$w^{[1]}$ and $w^{[2]}$ respectively. See Fig. \ref{fig:circle_rotation} where the two are compared.

When we compute with the first choice of $w$, we recover the truncated sum in the definition of the Birkhoff average. To estimate the error, we expect the difference $f(x_{N+1})-f(x_{N})$ to be of order one, implying that for $w_{equal}$, the error in the average is generally proportional to $N^{-1}$ in our figures.
The choice of a particular starting point also creates a similar uncertainty of order $1/N$. Every function $w$ is always positive  between $0$ and $1$. For all but the first choice, the function vanishes as $t$ approaches $0$ and $1$. In addition, going down the list, increasing number of derivatives of $w$ vanish for $t \to 0$ and $t\to 1$, with all derivatives of $w_{\exp}$ vanishing at $0$ and $1$. We thus expect the effect of the starting and endpoints to decay at the same rate as this number of vanishing derivatives. Indeed, we find that $w_{quad}$ corresponds approximately to order $1/N^2$ convergence, $w_{\sin^2}$ to $1/N^3$ convergence, and $w_{\exp}$ to convergence faster than any polynomial in $1/N$, i.e., for every integer $m$, there is a constant $C>0$ such that for $N$ sufficiently large, $|\Q_Nf-\int fd\mu|\leq CN^{-m}$. Figs.~\ref{fig:2DMap_overview}(b) and ~\ref{fig:mid_circle2} show this effect. We have not tried other $C^\infty$ weighting functions.

\subsection{Estimating error when the true rotation number is known}
In order to test the error in the calculation of rotation number, 
we present two examples below where we know the exact rotation number. This allows us to 
determine the actual error in the calculation for the $\Q_N$ method as $N$ increases. 
In both cases the error decreases to 
less than $10^{-31}$ and then it grows as $N$ increases, apparently due to accumulated 
round-off error. 

\begin{figure}
\centering
\subfigure[ ]{\includegraphics[width = .48\textwidth] {\Path 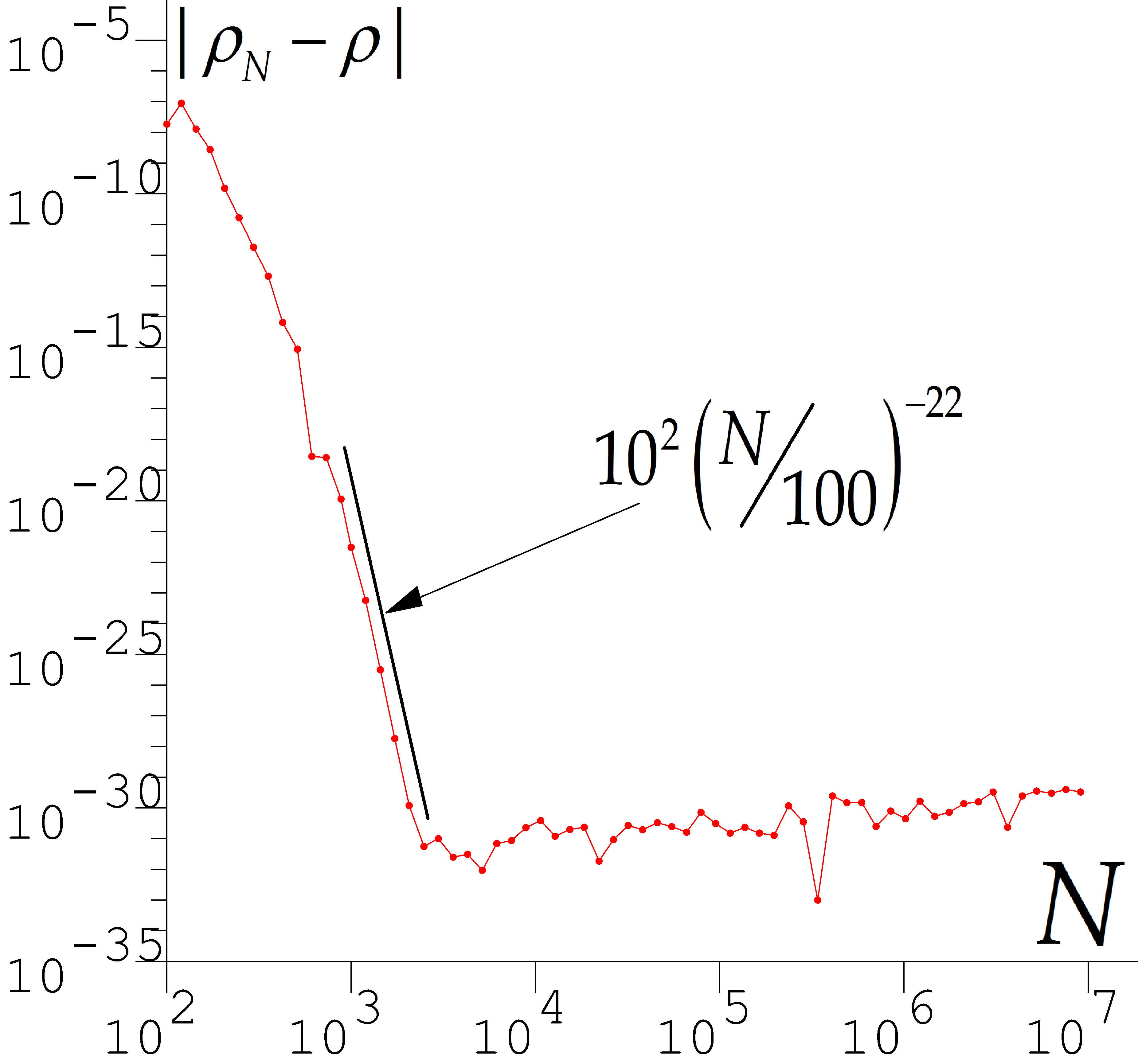}}
\subfigure[ ]{\includegraphics[width = .48\textwidth] {\Path 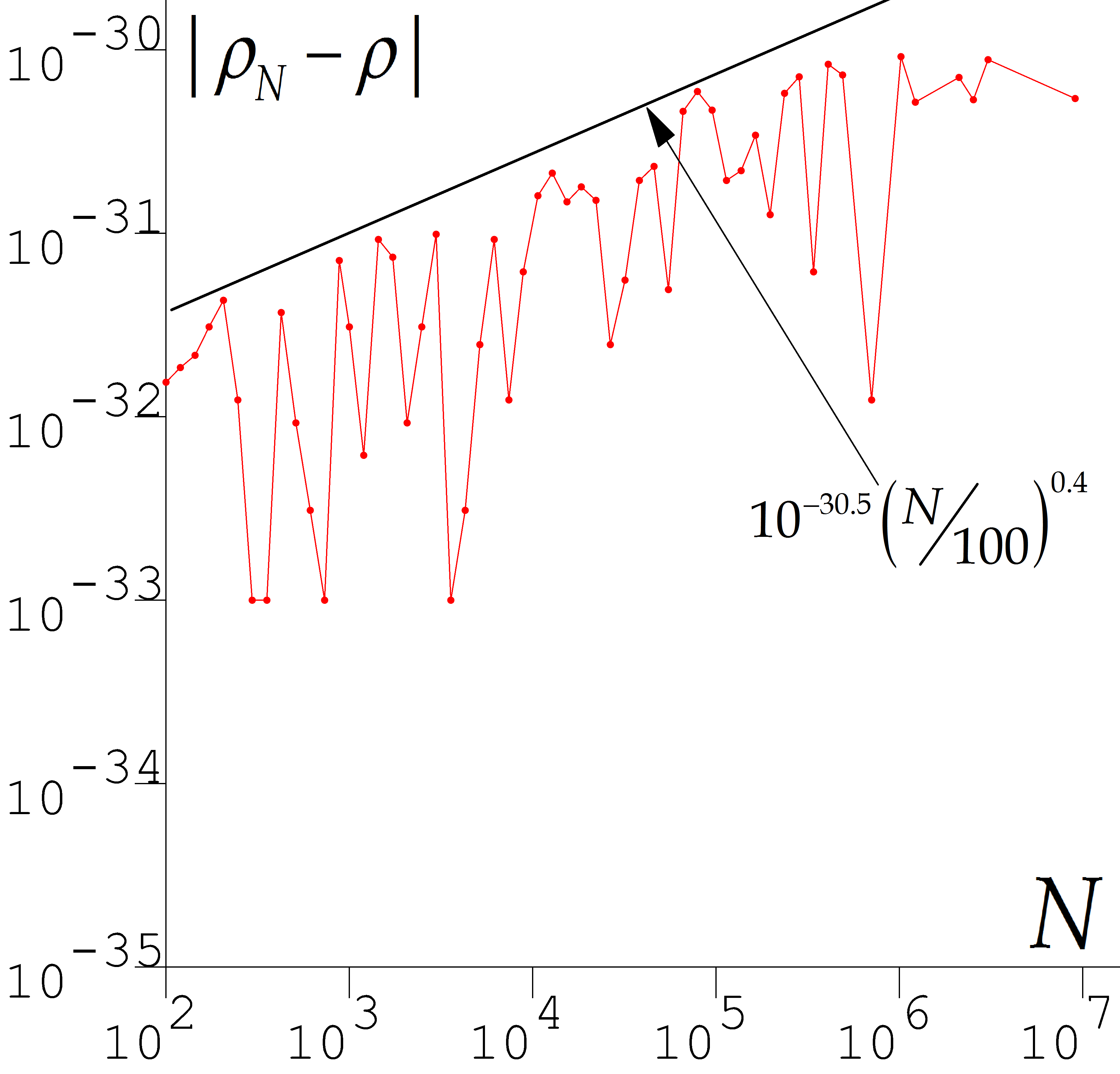}}
\caption{{\bf Testing how well the $\Q_N$ method can determine the rotation number.} Panel (a) shows the convergence in the calculation of a known rotation number $\rho=\sqrt{2}-1$, for the trajectory $(\phi_n)$ from Eq. (\ref{eqn:nonlinear_rotation}), with $\alpha=0.1$, $\beta=0.2$, and $\rho=\sqrt{2}-1$. The error quickly drops to the limit of numerical precision and then  increases slowly as $N$ increases. This increase in the error is apparently due to accumulated round-off error. Panel (b) shows the increasing round-off error in the rotation number for the trivial case ($\alpha=\beta=0)$. Here, $\phi_n=n\rho$ (mod $1$) and the error grows after attaining a minimum as $N$ increases. }
\label{fig:intrinsic_error}
\end{figure}

{\bf Example 1.} 
Let $(\theta_n)$ be an orbit under the rigid rotation described in Eq. \ref{eqn:Rigid_rotation} for a rotation by $\rho = \sqrt 2 -1$. Assume that what we observe is $\phi$, a perturbed version of $\theta$, namely, 
\begin{equation}\label{eqn:nonlinear_rotation}
\phi_{n}=\theta_n+\alpha\cos(2\pi\theta_n)+\beta\sin(2\pi\theta_n),\mbox{ where }\theta_n= n\rho\pmod 1.
\end{equation}
We use the Weighted Birkhoff average as in Eq. (\ref{eqn:weighted_rho})  (changing $y$ to $\phi$) to obtain an estimate of the rotation number $\rho$ from this orbit.
Fig. \ref{fig:intrinsic_error} shows the results for  $\alpha=0.1$ and $\beta= 0.2$ in (a) and for the case  $\alpha=0.0$ and $\beta= 0.0$ in (b). 

\begin{figure}
\centering
\subfigure[ ]{\includegraphics[width = .48\textwidth] {\Path 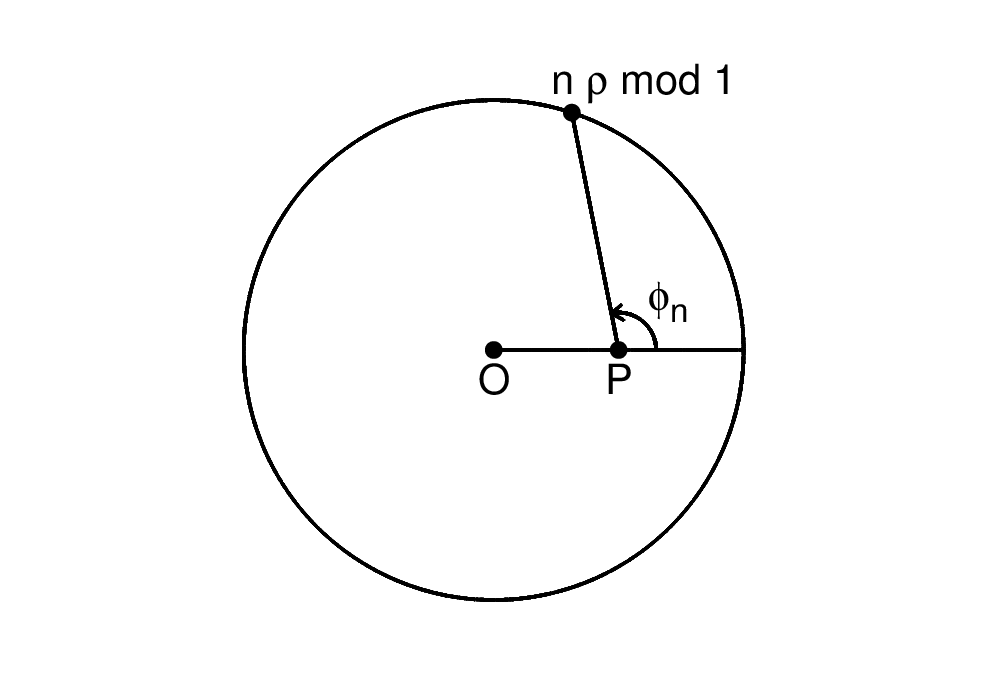}}
\subfigure[ ]{\includegraphics[width = .48\textwidth] {\Path 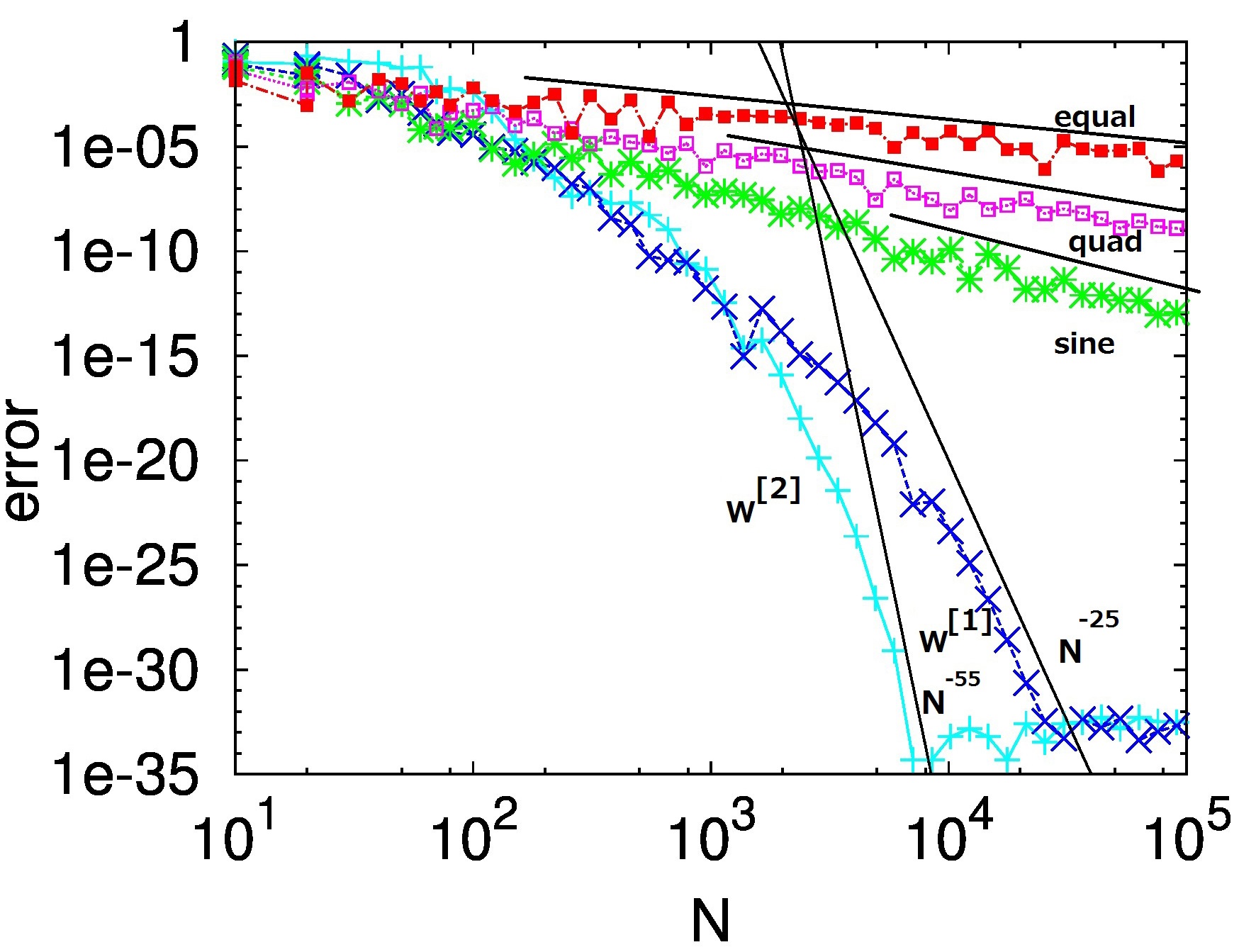}}
\caption{{\bf The error in the computed value of rotation number when the rotation number is known.} Panel (a) shows the geometric configuration of the problem, with a constant rotation vector $\rho = \sqrt 2 -1$ about the origin. The trajectory is $(n\rho \pmod 1)$, but the observer measures the angle $\phi$ as seen from its perspective at $P$, which is midway between the center of the circle, $O$, and the circle itself. 
Panel (b) shows the convergence in the rotation number calculation 
averaging $(\phi_{n+1}-\phi_n)$ five ways, using the Birkhoff average (top curve denoted ``equal'') and next ``quad'' and then ``sine''  and then the Weighted Birkhoff averages $\Q_N = \Q_N^{[1]}$ and $\Q_N^{[2]}$ (lowest curve), for the trajectory $(\phi_n)$ described in (a). The error from the known value $\sqrt 2 -1$ is calculated for several values of number $N$. The Weighted Birkhoff average $\WB_N^{[1]}$ reaches 32-digit accuracy by $N=30,000$ while $\WB_N^{[2]}$does so by $N=8,000$, at which point the $sin^2(\pi t)$ curve has an accuracy to $10^{-9}$, 
and its curve is proportional to $N^{-2.5}$.}
\label{fig:circle_rotation}
\end{figure}

{\bf Example 2.} Fig. \ref{fig:circle_rotation} shows a geometric version of the  problem from the previous example, and again the error in the rotation number is small. 

\subsection{Fourier coefficients and change of coordinates reconstruction}\label{sec:Fourier}

For a quasiperiodic curve as shown in Fig. \ref{fig:results_3B}(a), there are two approaches to representing the curve. Firstly, we can write the coordinates $(X,Y)$ as a function of $\theta\in S^1$, or secondly, we can reduce the dimension and represent the points on the curve by an angle $\phi\in S^1$, that is, $\phi(X(\theta),Y(\theta))$, which is also $V(\theta)=\theta+g(\theta)$. 
We have shown $g$ in Fig. \ref{fig:results_3B}(b) and the exponential decay of the norm of the Fourier coefficients in Fig. \ref{fig:results_3B}(d). 
We report the Fourier series for the periodic part $g(\theta)$.

Given a continuous periodic map $f:S^1\rightarrow\mathbb{R}$, the Fourier series representation of $f$ is the following:
\begin{equation}\label{eqn:sine_cosine_expansion}
\mbox{For every }\theta\in S^1,\ f(\theta)= \sum_{k=-\infty}^\infty a_k e^{i2\pi k\theta},
\end{equation}
where the complex coefficient $a_k$ is given by the formula
\begin{equation}\label{eqn:a_n}
a_k = \int_{\theta\in S^1} f(\theta) e^{-i2\pi k\theta} d\theta.
\end{equation}

If we only have access to an ergodic orbit $(x_n)$ on a curve, then we cannot use the fast Fourier transform 
as we only have the function values $f(x_n)$ along a quasiperiodic trajectory, and a rotation number $\rho$. Using interpolation to get the grid needed to apply a fast Fourier transform introduces significant interpolation errors. So instead, we obtain these coefficients using a Weighted Birkhoff average on a trajectory $(x_n)$ by applying the functional $\Q_{N}$. For $k=0$, we find $a_0$ by applying $\Q_N$ to the function $1$. Note that for all $k$, $a_{-k}=\bar{a_{k}}$. For $k>0$, we find $a_k$ as follows:
\begin{equation}\label{eqn:a_n_from_Q}
a_k = \Q_{N} (f(\theta) e^{-i2\pi k\theta}) = \sum_{n=0}^{N-1} f(x_n) e^{-i2\pi kn\rho} \hat{w}_{n,N}.
\end{equation}
This is depicted for the R3BP in Fig.~\ref{fig:results_3B}, for the standard map in Fig.~\ref{fig:StdMap}, for the forced van der Pol equation in Fig.~\ref{fig:vdP_global}. In all three one-dimensional cases, we depict $|a_k|$ as a function of $k$ for $k\geq 0$ only, as for all $k$, $|a_{-k}|=|a_k|$.  Our main observation is that the Fourier coefficients decay exponentially; that is, for some positive numbers $\alpha$ and $\beta$, in dimension one, the Fourier coefficients satisfy 
\begin{equation}\label{eqn:Fourier_decay}
|a_k| \le \alpha e^{- \beta |k|} \mbox{ for all }k\in\mathbb{Z}. 
\end{equation}
This is characteristic of analytic functions. 
Therefore, the conjugacy functions of all our examples are effectively, analytic,
``effectively'' meaning within the precision of our quadruple precision numerics.
In two dimensions, the computation of Fourier coefficients is similar, but instead of only having one exponential functions, for each $(j,k)$, we have two linearly independent sets of exponentials. 
\[ e^{i (j x + k y)} \mbox{ and } e^{i (j x - k y)}. \]
We define $a_{j,k}$ and $a_{j,-k}$ to be the complex-valued coefficients corresponding to these two functions. 

\subsection{Examples}
{\bf The standard map.} 
The standard map is an area preserving map on the two-dimensional torus, often studied as a typical
example of analytic twist maps (see \cite{StdMap1}). It is defined as follows
\begin{eqnarray}\label{eqn:StdMap}
S_1
\left(
\begin{array}{c}
x \\
y
\end{array}
\right)
=
\left(
\begin{array}{c}
x+y \\
y+  \alpha \sin (x+y)
\end{array}
\right)
\pmod{2\pi}.
\end{eqnarray}
\begin{figure}
\centering
\subfigure[ ]{\includegraphics[width = .35\textwidth] {\Path 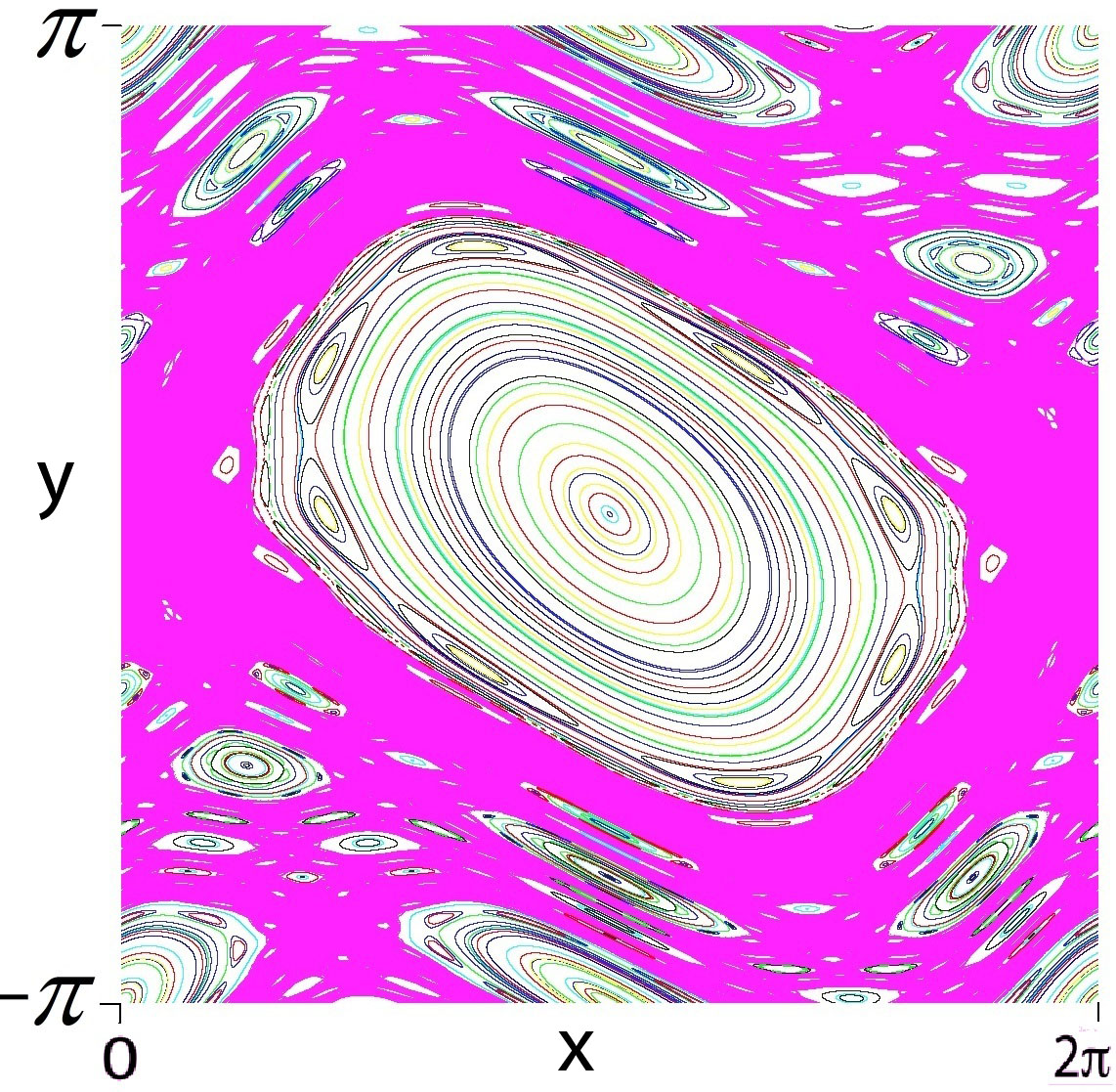}}
\subfigure[ ]{\includegraphics[width = .35\textwidth] {\Path 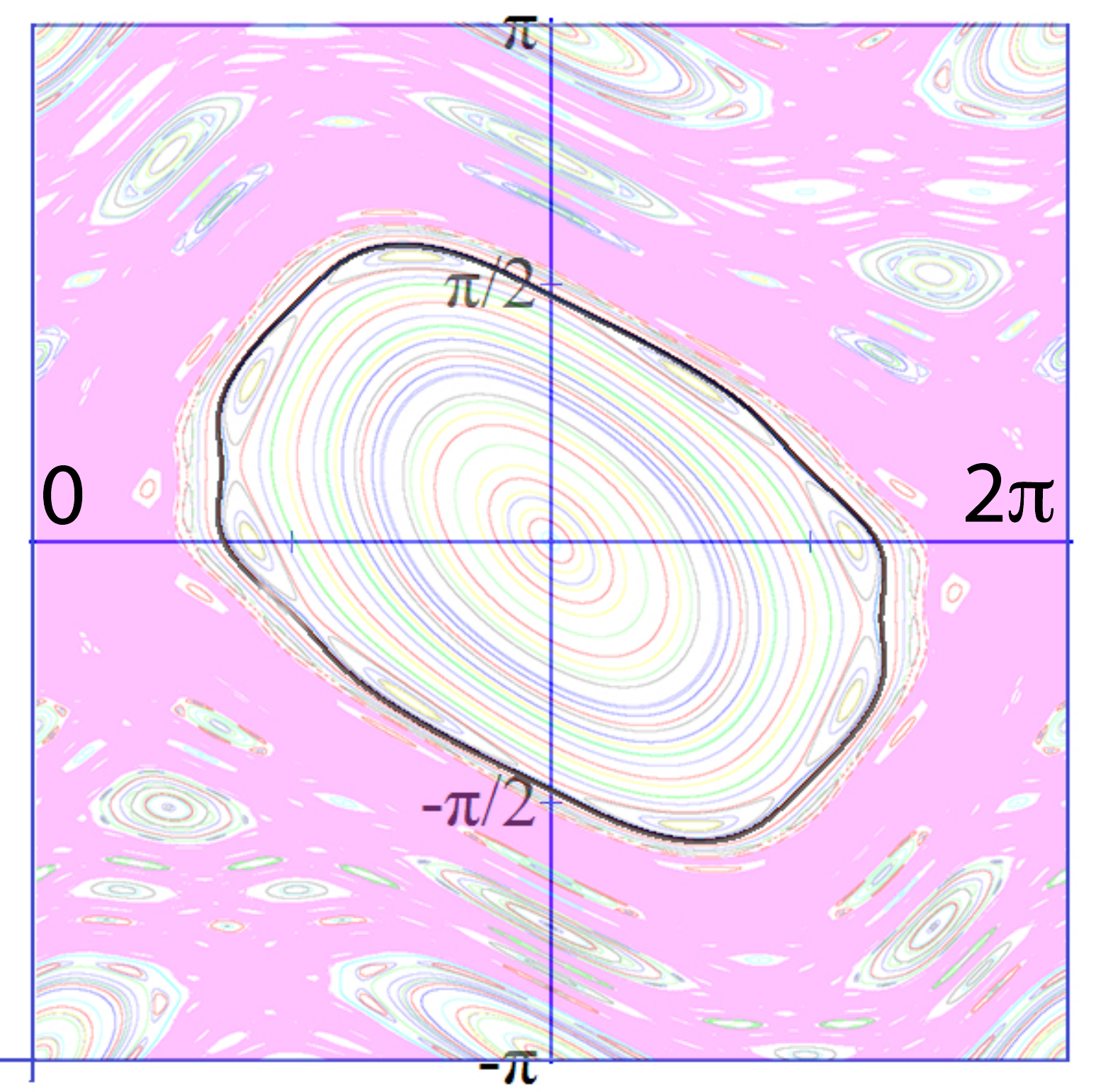}}
\caption{{\bf The standard map.} Panel (a) shows a variety of orbits from different initial conditions in the standard map $S_1$ defined in Eq.~\ref{eqn:StdMap}. We can see both chaos (shaded area) and quasiperiodic orbits under this map. A single curve with quasiperiodic behavior is plotted in panel (b). The orbit has initial conditions $(x,y) \approx (-0.607, 2.01)$. That is, if we restrict the map to this invariant curve, then it appears to be topologically conjugate to a rigid irrational rotation. }
\label{fig:StdMap_global}
\end{figure}
\begin{figure}
\centering
\includegraphics[height=.4\textwidth] {\Path 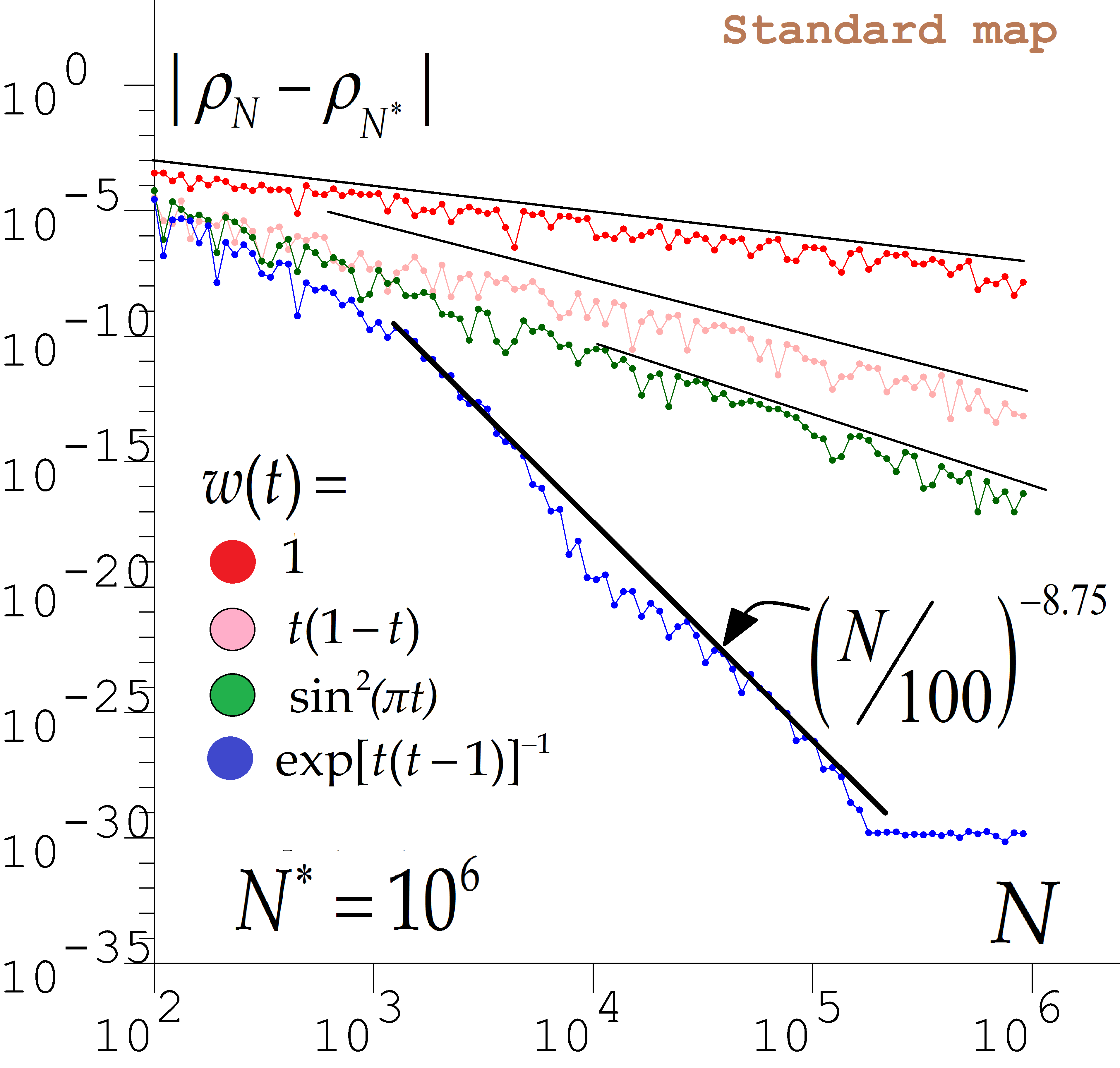}
\caption{\textbf{Rate of convergence in the rotation number for different weighting functions}. 
This shows results for the standard map using the same trajectory for each plot. For a given $w$ and a given number of iterates $N$, the rotation number $\hat{\rho}$ approximation is calculated for the curve  all using the same trajectory. The plotted error of the calculation is the difference $|\rho-\hat{\rho}|$ as a function of $N$. The exponential weighting function $w_{exp}$ allows $\Q_N$ to reach a limit by approximately $N=150,000$ at which point it is about $10^{16}$ times more accurate than the next best curve, the $\sin^2(\pi t)$ weighting. After that the error fluctuates by approximately $10^{-30}$. The black line indicates the average slope of the bottom-most curve for large $N$ and hence the approximate rate of
convergence for the $C^\infty$ Weighted Birkhoff average. The $\sin^2(\pi t)$ weighting curve convergence rate is proportional to $N^{-3.0}.$}
\label{fig:mid_circle2}
\end{figure}
\begin{figure}
\centering
\subfigure[ ]{\includegraphics[width = .35\textwidth] {\Path 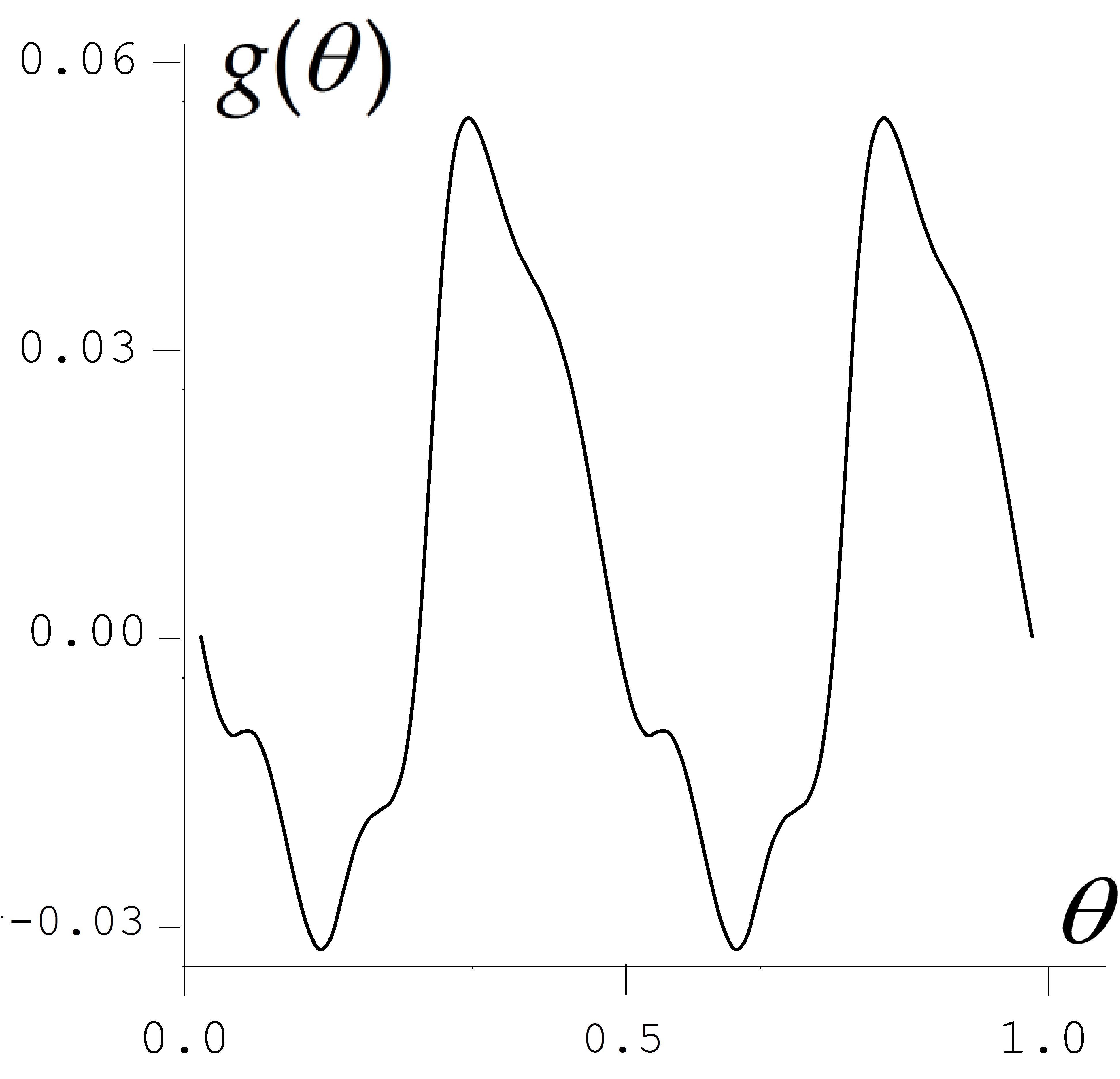}}
\subfigure[ ]{\includegraphics[width = .35\textwidth] {\Path 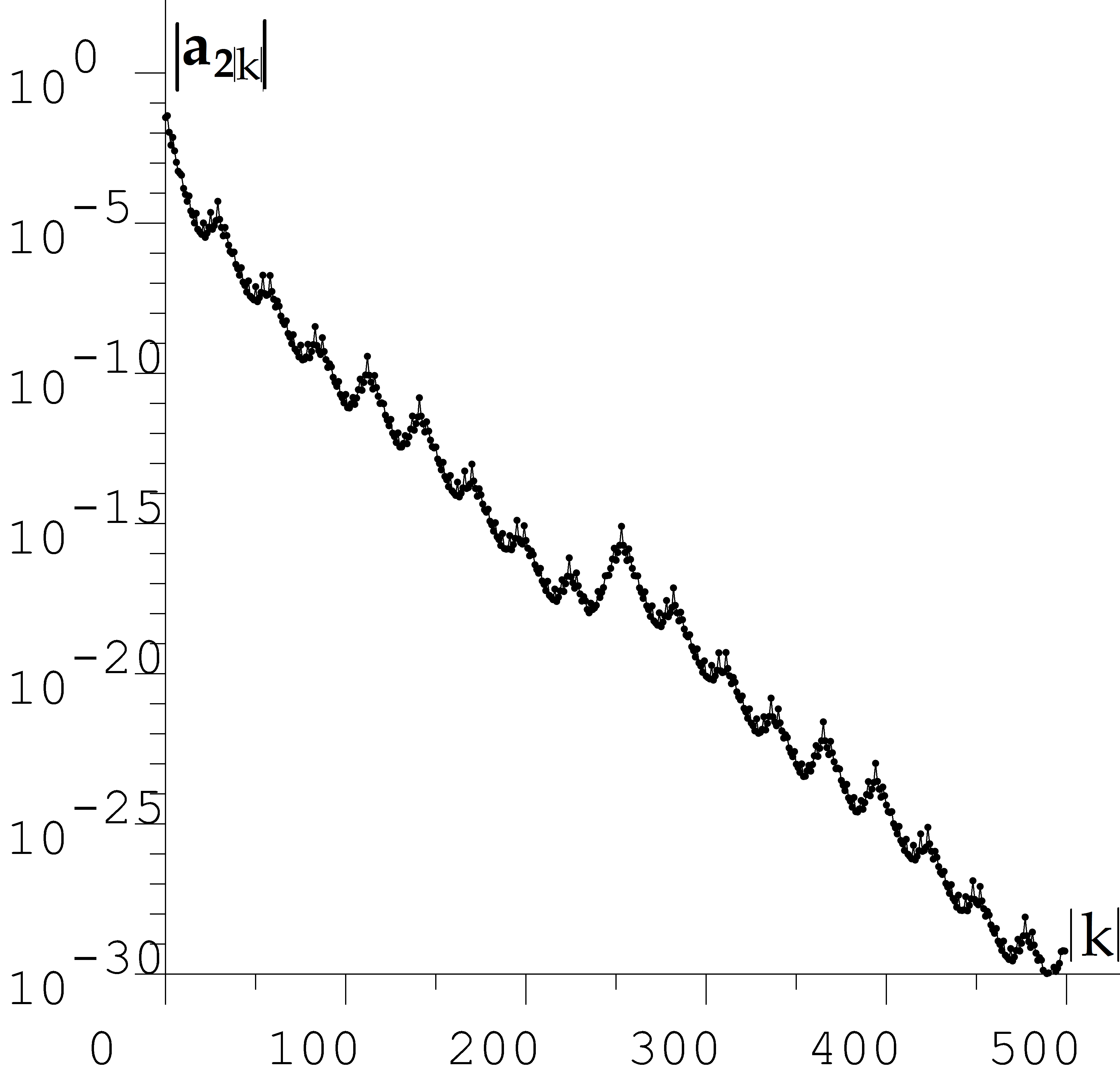}}
\caption{\textbf{The standard map conjugacy.} This figure shows the analysis of the quasiperiodic trajectory in Fig.~\ref{fig:StdMap_global}. Panel (a) depicts the periodic part $g(\theta)$ of the conjugacy between the quasiperiodic behavior and rigid rotation by $\rho$, measured with respect to the center point $(x,y)=(\pi,0)$. See Eq. \ref{eqn:g_conjugacy} for a description of $g(\theta)$. Panel (b) shows the decay of the Fourier coefficients. Since the conjugacy is an odd function, the odd-numbered Fourier coefficients are zero and therefore have been omitted from the picture. The decay of the Fourier terms can be bounded from above be an exponential decay, which suggests that the conjugacy is analytic. An orbit of length $N=10^{7}$ is used for these computations. A smaller orbit of length $N=10^6$ does not lead to any significant changes.} 
\label{fig:StdMap}
\end{figure}
In this paper, we only consider the case $\alpha=1.0$. Fig.~\ref{fig:StdMap_global}(a) shows the trajectories starting at a variety of different initial conditions plotted in different colors. The shaded set is a large invariant chaotic set with chaotic behavior, but many other invariant sets consist of
one or more topological circles, on which the system has quasiperiodic behavior. For example, initial condition $(\pi,1.65)$ yields chaos while $(\pi,1.5)$ yields a quasiperiodic trajectory. As is clearly the case here, one-dimensional quasiperiodic sets often occur in families for non-linear processes, structured like the rings of an onion. There are typically narrow bands of chaos between quasiperiodic onion rings. Usually these rings are differentiable images of the $d$-torus. We have computed the rotation number to be $0.12055272197375513300298164369839$ for one such  standard map orbit shown in the Fig.~\ref{fig:StdMap_global}(b) using quadruple precision. 
 Fig. \ref{fig:mid_circle2} shows a convergence rate of $O(N^{-8.75})$ in computing the rotation number using $\WB_N$.
Fig.~\ref{fig:StdMap}(a) shows the periodic parts $g$ of conjugacies of the  quasiperiodic orbit, and Fig.~\ref{fig:StdMap}(b) shows the absolute values of Fourier coefficients representing an exponential decay.

\begin{figure}
\centering
\subfigure[ ]{\includegraphics[width=.3\textwidth] {\Path 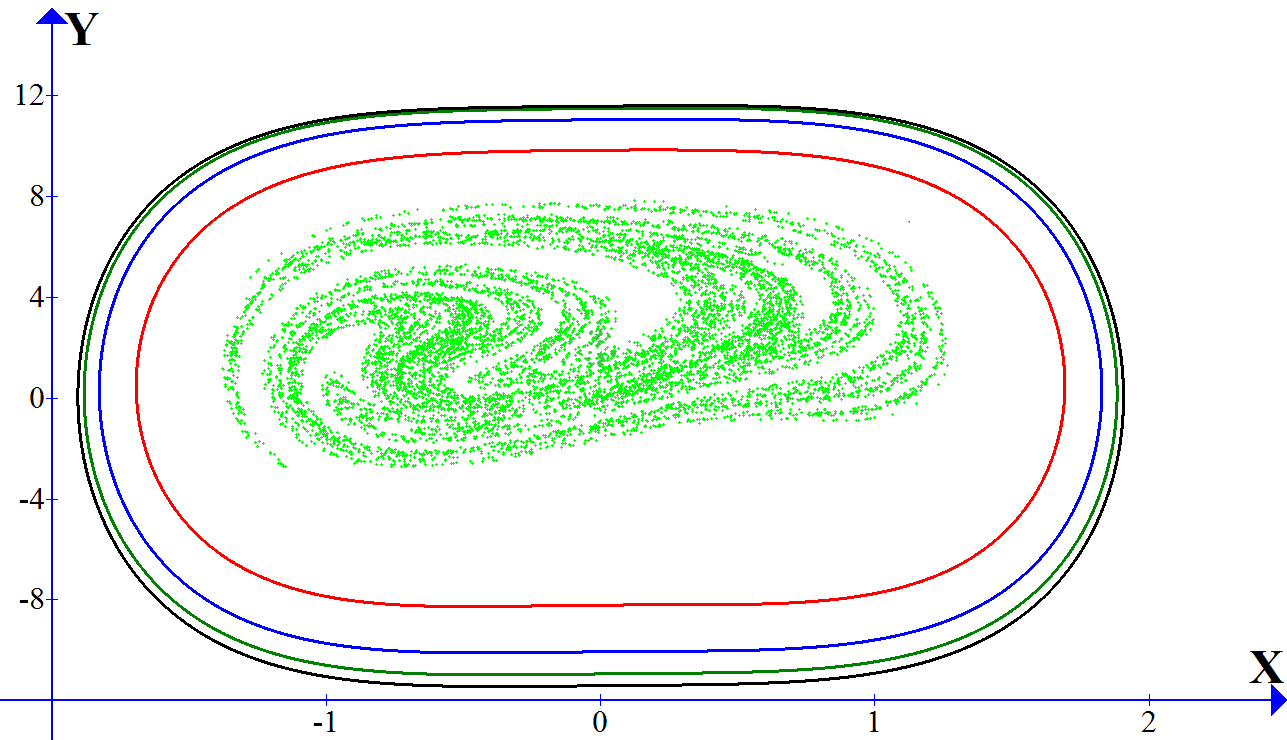}}
\subfigure[ ]{\includegraphics[width = .3\textwidth] {\Path 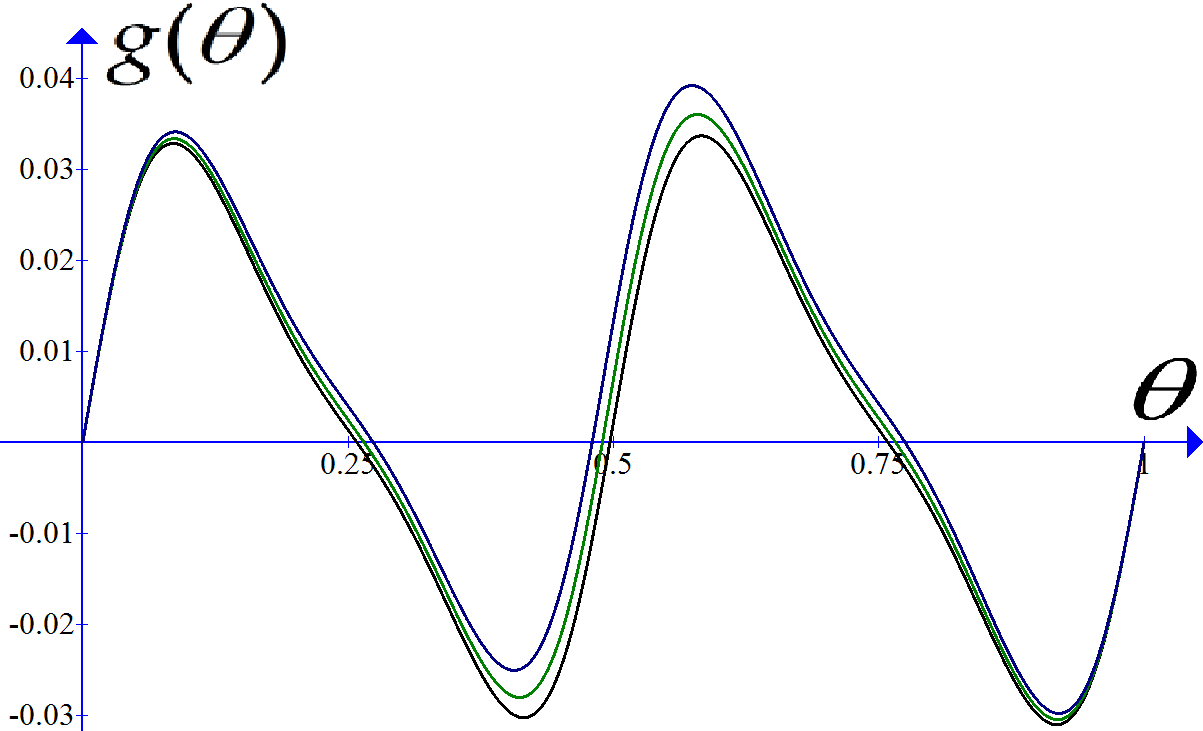}}
\subfigure[ ]{\includegraphics[width = .3\textwidth] {\Path 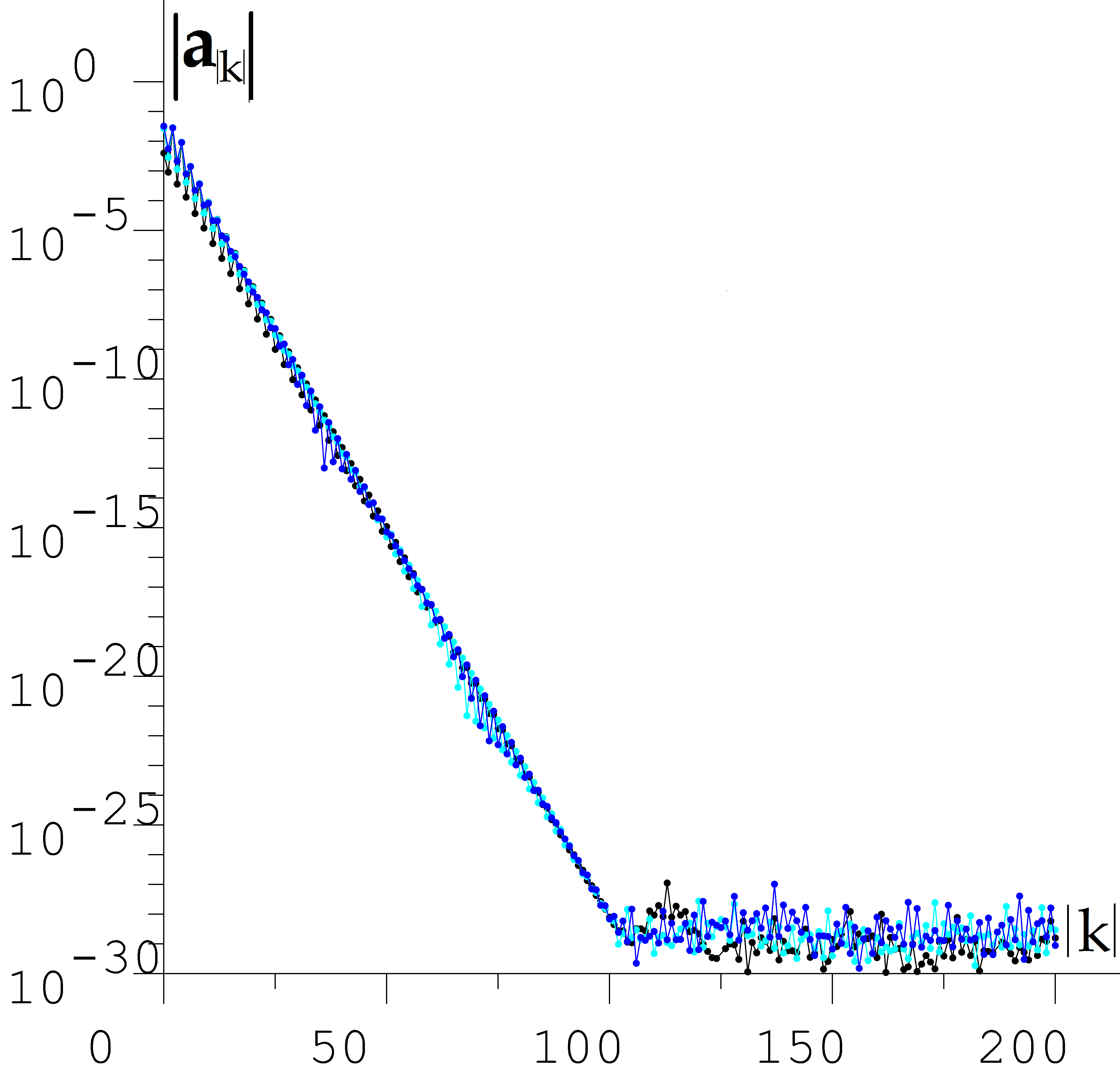}}
\caption{ {\bf Forced van der Pol oscillator.} Panel (a) shows attracting orbits for a number of different forcing values $F$ for the stroboscopic map of the van der Pol flow given in Eq. \ref{vdPol_forced}. The plot depicts points $(X,Y) =(x(t_k),x'(t_k))$, where $t_k = 2 k \pi/0.83, k = 0,1,2,\dots.$ The chaotic orbit lying inside the cycles corresponds to $F=45.0$. There are stable quasiperiodic orbits shown as curves, which from outermost to innermost correspond to $F=5.0$, $15.0$, $25.0$ and $35.0$ respectively. Panel (b) is the periodic part $g(\theta)$ of the conjugacy (Eq.~\ref{eqn:conjugacy1D}) to a rigid rotation, measured with respect to the center point $(0,0)$ for $F= 5.0, 15.0, $and $25.0$. Panel (c) shows the absolute values of the Fourier coefficients of $g$. Their linear decay as $|k|$ increases in this log-linear plot means the coefficients decay exponentially fast, which is a signature of the analyticity $g$ and hence of the conjugacy to a rigid rotation. The decay is exponentially fast down to the resolution of the numerics.}
\label{fig:vdP_global}
\end{figure}

{\bf The forced Van der Pol oscillator.} Fig.~\ref{fig:vdP_global}(a) shows attracting orbits for the time-$2\pi/0.83$ map of the following periodically forced Van der Pol oscillator with nonlinear damping~\cite{vdPol_forced}
\begin{equation}\label{vdPol_forced}
\frac{d^2 x}{dt^2}-0.2 \left( 1-x^2 \right) \frac{dx}{dt}+20 x^3=F\sin \left(0.83 t \right),
\end{equation}
for several values of $F$. While the innermost orbit shown is chaotic, the outer orbits are topological circles with quasiperiodic behavior\footnote{As with the standard map, we have specified all non-essential parameters rather than stating the most general form of the Van der Pol equation. Our computed rotation numbers for the three orbits $F=5.0$, $15.0$, and $25.0$ are $0.29206126329199589285577578718959$, $0.37553441113144010884908928083318$ and $0.56235370092685056634419221336154$ respectively.}.
Each curve was assigned the angular coordinates of Eq. \ref{eqn:rot_num_defn} by assigning the points on the curve the angle with respect to the origin $(0,0)$.
Fig.~\ref{fig:vdP_global}(b) shows the periodic parts $g$ of conjugacies of quasiperiodic orbits, and Fig.~\ref{fig:vdP_global}(c) shows the absolute values of Fourier coefficients representing exponential decays.

{\bf A two-dimensional torus map.} So far, the quasiperiodic sets studied here are closed curves. We now describe an example \cite{Grebogi:83,Grebogi:85,Kim} of a two-dimensional quasiperiodic torus map on $\mathbb{T}^2$.  This is a two-dimensional version of Arnold's family of one-dimensional maps (see \cite{Arnold}).  The map is given by $(T_1,T_2)$ where 
\begin{eqnarray}
T_1(x,y) = \left[ x+ \omega_1 + \frac{\epsilon}{2 \pi} P_1(x,y) \right] \pmod 1, \label{eqn:2D_map} \\
T_2(x,y) = \left[ y+ \omega_2 + \frac{\epsilon}{2 \pi} P_2(x,y) \right] \pmod 1,\nonumber
\end{eqnarray}
and $P_i(x,y), i=1,2$ are periodic functions with period one in both variables, defined by:
\[
P_i(x,y) = \sum_{j=1}^4 a_{i,j} \sin (2 \pi \alpha_{i,j}), \mbox{ with } \alpha_{i,j}= r_j x + s_j y +b_{i,j}.
\]
The values of all coefficients are given in Table~\ref{table}. 
\begin{table}
\begin{center}
{\begin{tabular}{ |c | c| }
\hline
Coefficient & Value \\
\hline
$\epsilon$ & $ 0.4234823$ \\[1pt]
$\omega_1$ & $ 0.71151134457776362264681206697006238 $ \\[2pt]
$\omega_2$ & $ 0.87735009811261456100917086672849971 $ \\[1pt]
$a_{1,j}$ & $(-0.268,-0.9106,0.3,-0.04)$ \\[1pt]
$a_{2,j}$ & $(0.08,-0.56,0.947,-0.4003)$ \\[1pt]
$b_{1,j}$ & $(0.985,0.504,0.947,0.2334)$ \\[1pt]
$b_{2,j}$ & $(0.99,0.33,0.29,0.155)$ \\[1pt]
$r_j$ & $(1,0,1,0)$ \\[1pt]
$s_j$ & $(0,1,1,-1)$ \\[1pt]
Computed $\rho_1$ & $0.718053759982066107095244936117$ \\[1pt]
Computed $\rho_2$ & $0.885304666596099792113366824157$ \\[1pt]
\hline
\end{tabular}}
\end{center}
\caption{\label{table}
{\bf Coefficients for the torus map.} All values are used in quadruple precision, 
but in this table the repeated zeros on the end of the number are suppressed.}
\end{table}
This choice of this function is based on~\cite{Grebogi:83,Grebogi:85}. The papers use the same form of equation, though the constants are close to but not precisely the same as the ones used previously. This fits with the point of view advocated by these papers: that the constants should be randomly chosen. Since we are using higher precision, we have chosen constants that are irrational to the level of our precision. The forward orbit is dense on the torus, and the map is a nonlinear map which exhibits two-dimensional quasiperiodic behavior.

\begin{figure}
\centering
\subfigure[ ]{\includegraphics[height= .3\textwidth] {\Path 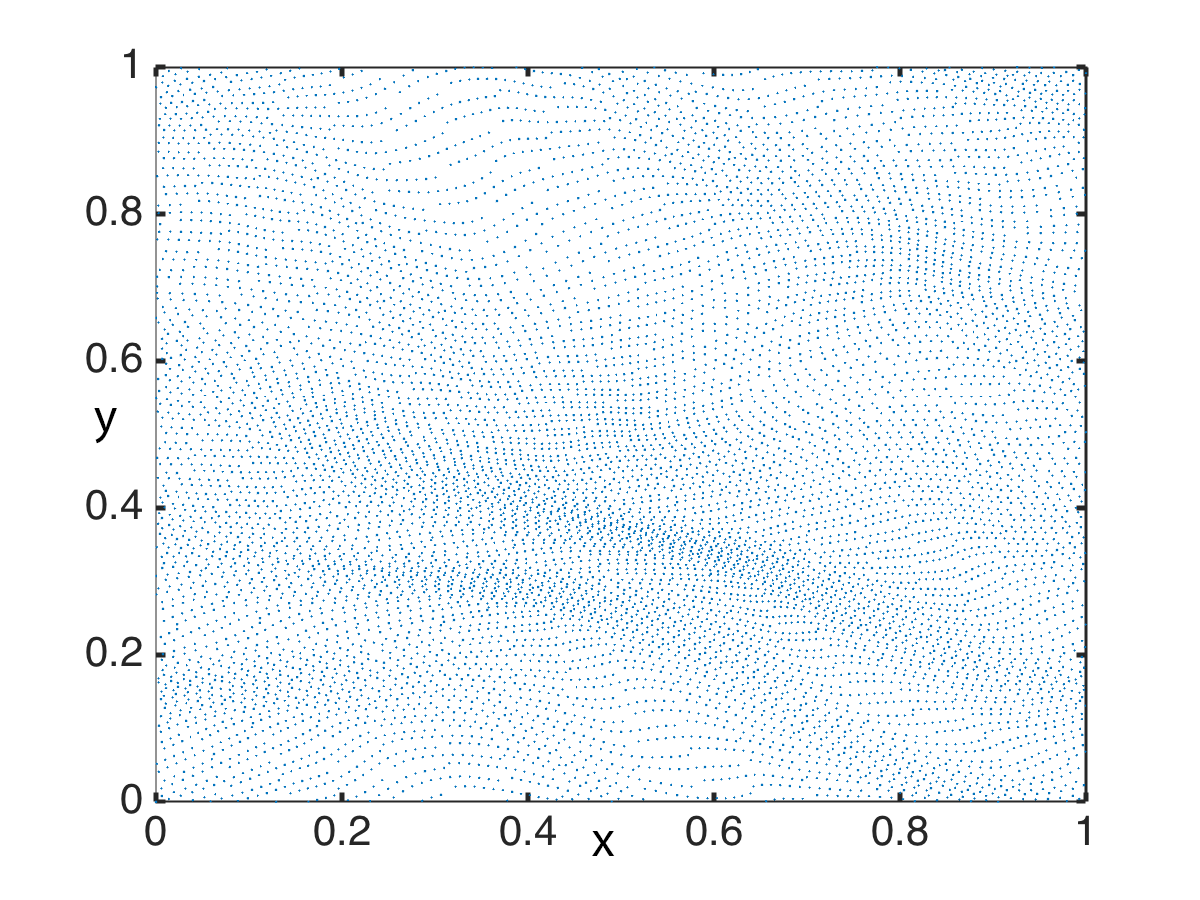}}
\subfigure[ ]{\includegraphics[height= .3\textwidth] {\Path 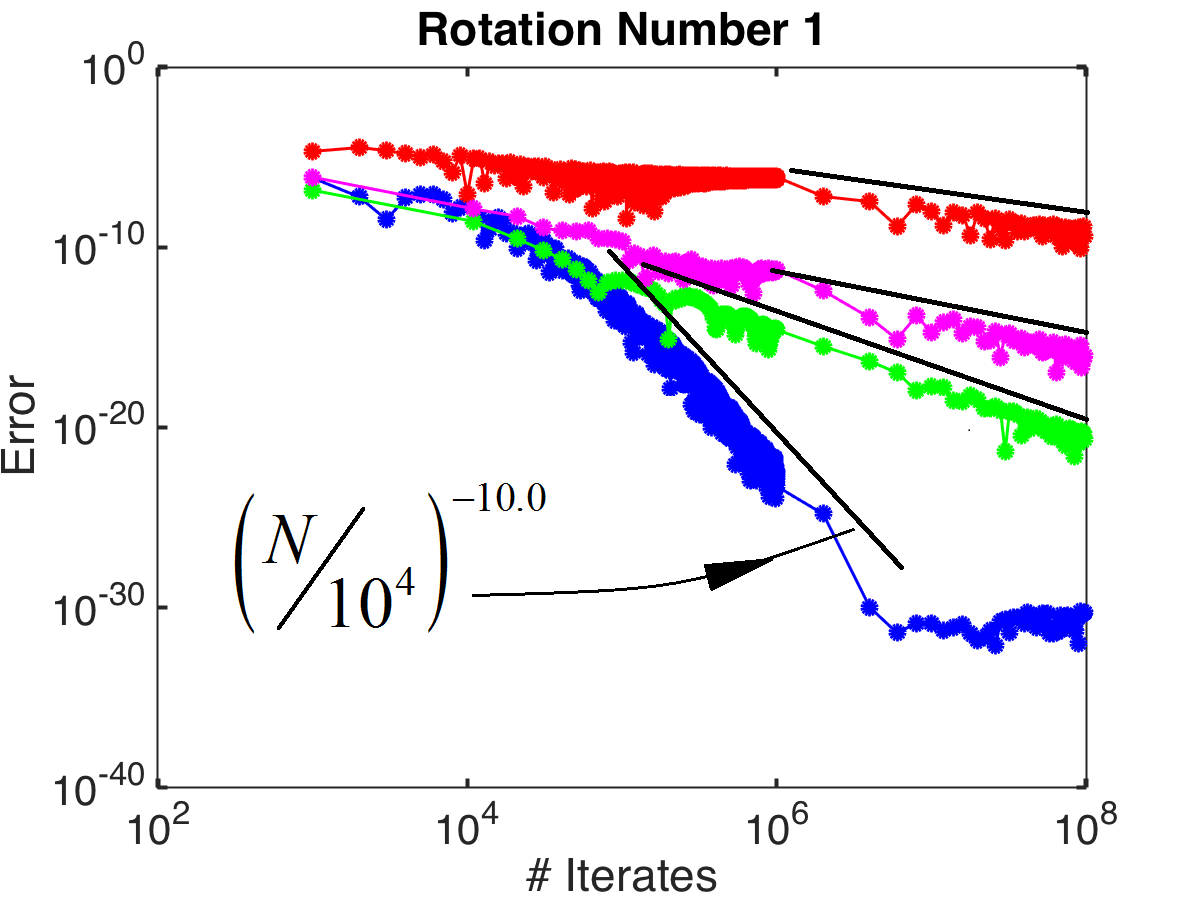}}
\caption{\textbf{Two-dimensional torus map.} Panel (a) shows an orbit of length $10^{4}$ for the two-dimensional quasiperiodic torus map. The orbit appears to be dense, which is consistent with quasiperiodicity. 
Panel (b) shows the convergence rate for the first rotation number for the four different weighting functions given in Eq. \ref{eqn:weights}, resp. from highest curve to lowest, plotted in red, magenta, green, and blue (online).
The black line is an upper bound for the blue curve and its exponent indicates the rate of convergence when using the $C^\infty$ weight. It has a slope of $-10$ in the log-linear graph while the corresponding slope for the $\sin^2(\pi t)$ weighting function would be $-3.$ That is, its rotation rate's error is proportional to $N^{-3}$.}
\label{fig:2DMap_overview}
\end{figure}
\begin{figure}
\centering
\subfigure[ ]{\includegraphics[width = .46\textwidth] {\Path 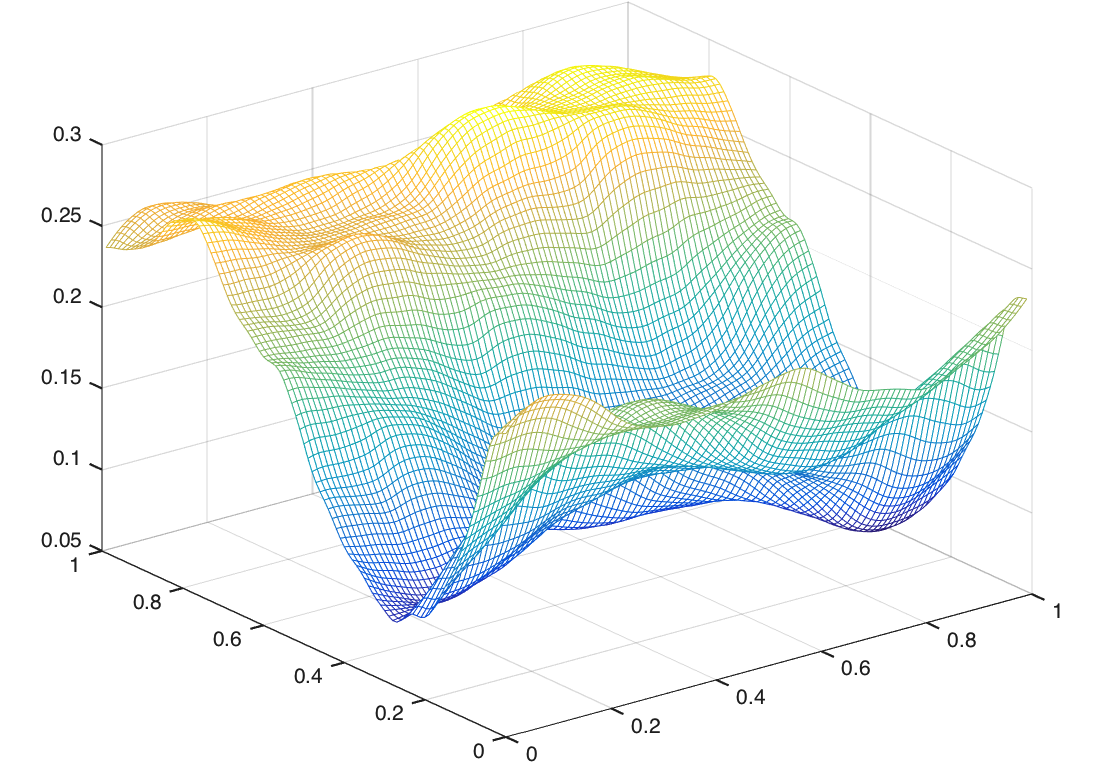}}
\subfigure[ ]{\includegraphics[width = .46\textwidth] {\Path 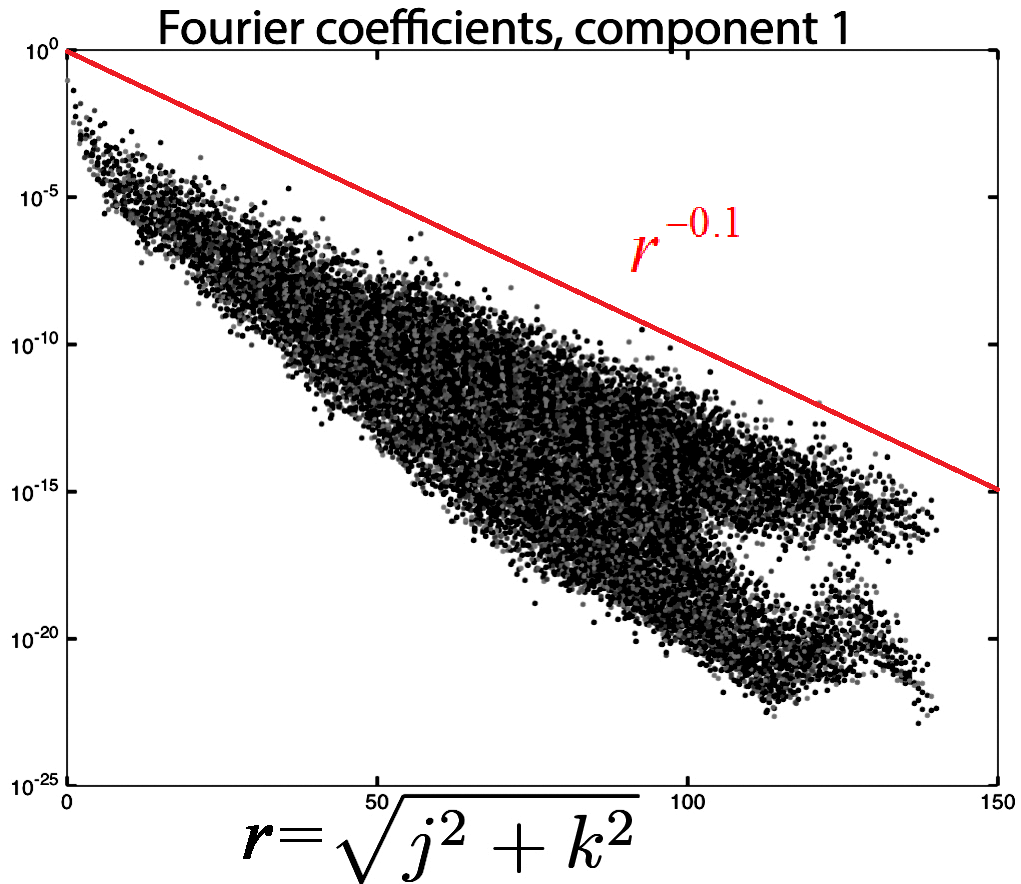}}
\caption{\textbf{Conjugacy for the two-dimensional torus map.} Panel (a) depicts the reconstruction of the periodic part $g$ (see Eq.  \ref{eqn:g_conjugacy}) of the first component of the conjugacy function for the torus map. The surface is colored by height. The surface is created using the Fourier coefficients shown in (b). The second conjugacy function is similar but not depicted here. Panel (b) shows the norm of the Fourier coefficients $|a_{j,k}|$ and $|a_{j,-k}|$, for $0 \le j,k \le 100$ for the first component of the conjugacy function from panel (a). }
\label{2DTorus_results}
\end{figure}

Fig.~\ref{fig:2DMap_overview}(a) depicts iterates of the orbit, indicating that it is dense in the torus. We use our Weighted Birkhoff average to compute the two Lyapunov exponents, which have super convergence to  zero. Fig.~\ref{fig:2DMap_overview}(b) shows one of them. In terms of method, this is just a matter of changing the function $f$ used in $\Q_N$ in Eq. \ref{eqn:QN}. 
Likewise, finding rotation numbers in two dimensions uses the same technique as in the one-dimensional case ({\em cf.} Fig.~\ref{fig:2DMap_overview}(c)). In all of our calculations, the computation is significantly longer than in dimension one in order to get the same accuracy, perhaps because in dimension two sufficient coverage by a trajectory may vary like the square of the side length of the domain. Fig. \ref{fig:2DMap_overview}(b) shows a convergence rate of $O(N^{-10.0})$.

The reconstructed conjugacy function for the two-dimensional torus is depicted in Fig.~\ref{2DTorus_results}(a). The decay of Fourier coefficients shown in Fig.~\ref{2DTorus_results}(b) displays $\sqrt{j^2 +k^2}$ on the horizontal axis, and $|a_{j,k}|$ on the vertical axis, where both of these coefficients are complex, meaning that $|\cdot|$ represents the modulus. Again here, the coefficients decay exponentially, though the decay of coefficients is considerably slower in two dimensions. 
The data set looks quite a lot more crowded in this case, since there are 
$10^4$ coefficients and many different values of $(j,k)$ have almost same values of $\sqrt{j^2 + k^2}$.
In addition, the coefficients $a_{j,k}$ generally converge at different exponential rates. This is why there is a strange looking set consisting of an upper and a lower cloud of data in Fig.~\ref{2DTorus_results}(b). While more information on the difference between these coefficients is gained by interactively viewing the data in three dimensions, we have not been able to find a satisfactory static flat projection of this data. We feel that in a still image, the data cloud shown conveys the maximum information.

{\bf The three-body problem revisited.} 
We have computed trajectories for the Poincar{\'e} return map using an $8^{th}$ order Runge-Kutta method with time step $2\times 10^{-5}$, in quadruple precision and  Fig.~\ref{fig:results_3B}(c) is consistent with 30-digit accuracy of the rotation number (See item 1 below.)
for the quasiperiodic orbit of the three-body problem labeled $B_1$ in Fig.~\ref{fig:3B1}(b). 
The initial condition for this orbit is $(q_1,q_2,p_1,p_2) = (-0.15, 0, 0, p_2)$ where $p_2\approx 5.41$ is chosen so that the Hamiltonian is about $-2.63$. 
While it is more straightforward to obtain a numerical trajectory when the quasiperiodic trajectory is stable, it is also possible when it is a saddle or a repeller. Our Weighted Birkhoff average $\Q_N$ approach works equally well for both cases. The extent of convergence of $\Q_N$ is limited by the accuracy of the trajectory data.We now list results of our numerical methods for the restricted three-body problem. 
\begin{enumerate}
\item The rotation number is $0.063961728757453097164077724400302$, computed to 30-digit accuracy, and 
Fig.~\ref{fig:results_3B}(c) shows the accuracy plateauing at about 30-digit accuracy.
\item We 
compute 200 terms of the Fourier series, the last 125 of which have magnitude near 0 (i.e., less than $10^{-30}$).
There is a conjugacy map $h$ between the 
first return map and a rigid rotation on the circle. Evaluating the Fourier series allows us to reconstruct the conjugacy map ({\em cf.} Fig.~\ref{fig:results_3B}(b)). 
\item
We find that the coefficients for the conjugacy map decrease exponentially fast (see Fig. \ref{fig:results_3B}(d)).
That is a signature of a real analytic function as in  Fig.~\ref{fig:results_3B}(b).
\item
The high rate of convergence of $O(N^{-15.0})$ for the rotation number in Fig. \ref{fig:results_3B}(c) suggests that we have an effective computational method that yields an accuracy that is close to the limit of numeric precision, provided $N$ is sufficiently large. 
\end{enumerate}

{\bf Speed of convergence of Fourier series for a conjugacy.} In a separate report \cite{Siegel}, we examine conjugacies of the quasiperiodic curves in the Siegel disk. The map is a simple one-dimensional complex dynamical system $z_{n+1}=f(z_n)$, where 
\begin{equation}
  f(z)=z^2+e^{2\pi i \rho}z  
\end{equation}
and $\rho=(\sqrt{5}-1)/2$ \cite{luque:Rot,Llave-Petrov2008}.


It is found (see Fig. 3 in \cite{Siegel}) that the invariant curves get more irregular near the boundary of the disk. While typical smooth quasiperioidc curves require about 70 Fourier coefficients for 30-digit precision, near the boundary of the Siegel disk the curves are much more irregular and can require 24,000 coefficients (or more) for the same precision. The curves thus become increasingly fractal looking near the boundary and the Fourier series converges much more slowly. We remark that we would expect a similar slower convergence when exploring a case  like the well known last KAM circle of the Standard Map.

{\bf Sources of error.}
We end this subsection by noting a few sources of error in the computation of Fourier coefficients. If the number of iterates $N$ is too small, then we will not have sufficient coverage to get a good approximation of the coefficients, and the problem becomes more acute as the coefficient number $|k|$ grows. 
If the approximation of the rotation number is not accurate, then we cannot expect the approximations of our Fourier coefficients to be good either, and given an error in the rotation number, there will be a $k_{\max}$ such that the Fourier coefficients
$a_k$ with $|k|>k_{\max}$ cannot be approximated with any reasonable accuracy. A more subtle form of a error comes from the fact that if the rotation number we are trying to estimate is close to being commensurate with the rotation number, then we will get unexpectedly insufficient coverage of the space when performing iteration. See Section \ref{sec:almost_rational}.
%

\subsection{Lyapunov exponents computed as a Weighted  average}\label{sec:Lyap}

Lyapunov exponents are an important characterization of the dynamics resulting from the map $T$. 
They measure the rate at which nearby trajectories diverge or converge and can be used to distinguish between chaos and quasiperiodicity, for example, the existence of a positive Lyapunov exponent implies chaos, while for quasiperiodic systems, all $d$ Lyapunov exponents for Eq. \ref{quasi} are zero. 
In this section, we will show a way to obtain ``super-convergence'' to the Lyapunov exponents of a quasiperiodic dynamics on a 2D torus. 
Since we are considering the dynamics restricted to the torus, we do not calculate the Lyapunov exponents for the normal sub-bundle.

\begin{figure}
\centering
\subfigure[ ]{\includegraphics[height= .38\textwidth]{\Path 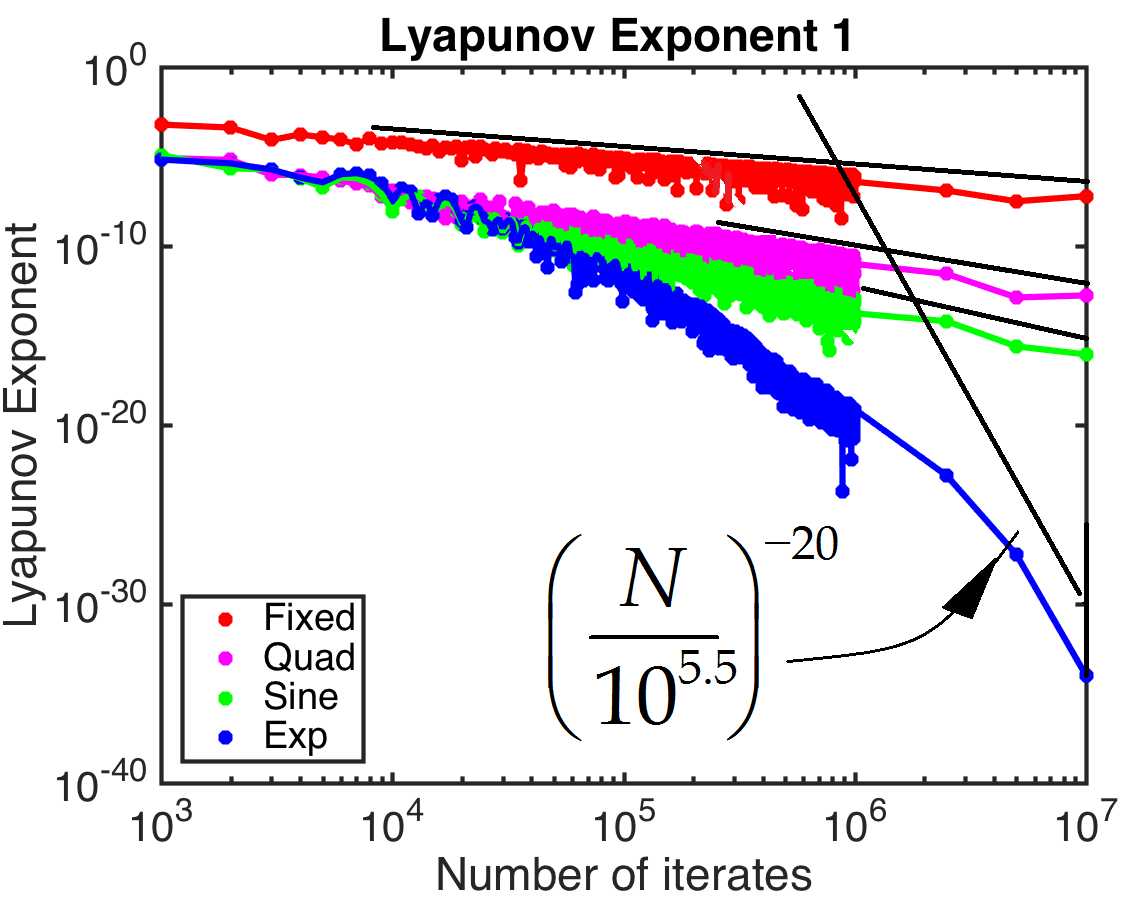}}
\subfigure[ ]{\includegraphics[height= .38\textwidth]{\Path 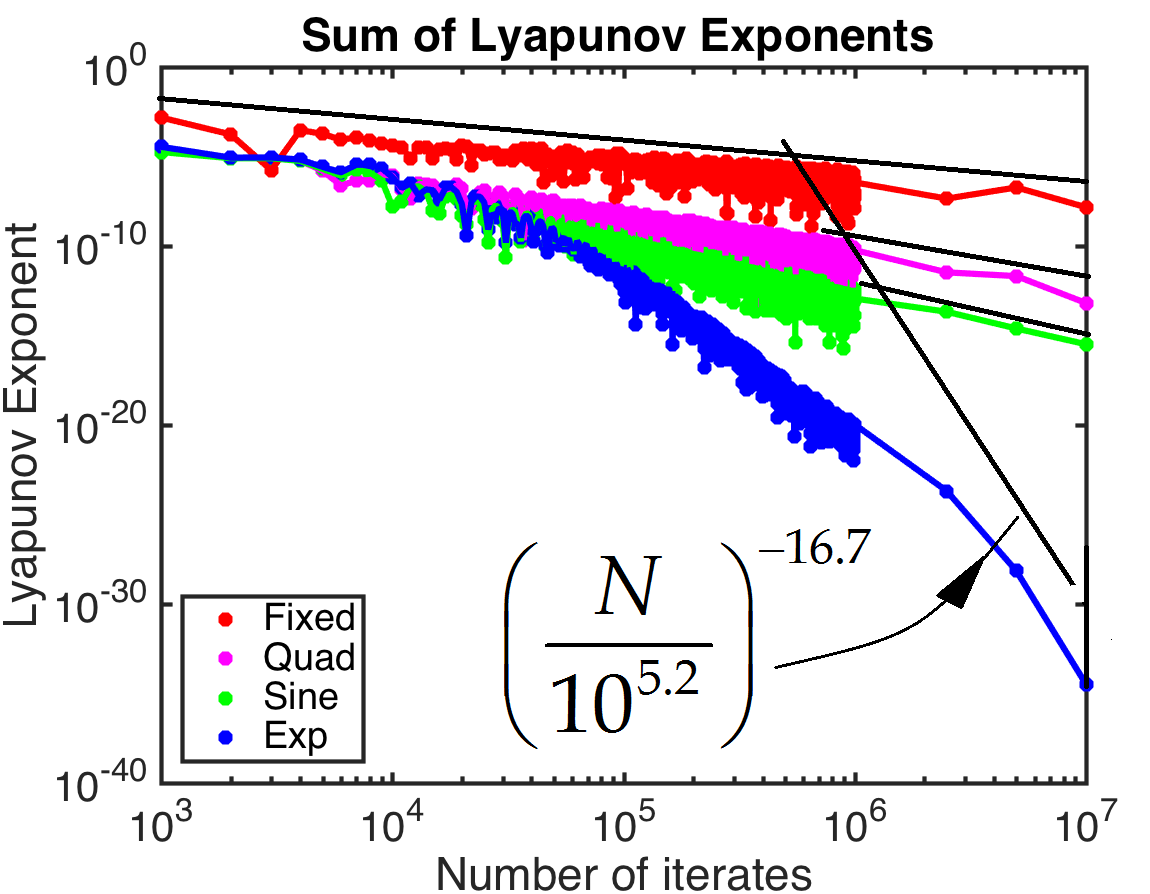}}
\caption{\textbf{Lyapunov exponents computed with $\WB_N$.}
The two Lyapunov exponents of the torus map from Fig.~\ref{fig:2DMap_overview} are computed using $\WB_N$ and we find them to be $0$ up to our numerical accuracy. (a) The computation for the first Lyapunov exponent is shown in blue for up to $10^7$ iterates. The other three curves depict the values attained using the three other weighting functions, using the  same order (and color scheme online) as in Fig.~\ref{fig:2DMap_overview}. Panel (b) shows the sum of the two Lyapunov exponents  computed via direct computation of the determinant as described in Eq. \ref{eqlyapsum}. In the two panels, the $\sin^2(\pi t)$ graphs are approximately proportional to  $N^{-3}$.}
\label{fig:2DMap_lyap}
\end{figure}
Lyapunov exponents are usually calculated numerically as an average  of logarithms  with all terms weighted equally, and the computations are therefore limited by extremely slow convergence rates. 
For quasiperiodic systems we would expect the error in the exponents to be $O(1/N)$. For 2-dimensional quasiperiodic systems, the Lyapunov exponent can be expressed as a Weighted Birkhoff average (using $\Q_N$) and we observe that we obtain the Lyapunov
exponents much faster. Figs.~\ref{fig:2DMap_lyap}(a) and \ref{fig:2DMap_lyap}(b) show convergence rates of $O(N^{-20.0})$ and $O(N^{-16.7})$ for the computation of the Lyapunov exponents using $\WB_N$. (We make no claim that we can prove convergence in (a) is super fast.) Using the
Weighted Birkhoff average does not change the limit.  

Recall that given a $d$-manifold $M$ and a map $T:M\to M$ with an invariant probability measure $\mu$, Oseledets' multiplicative ergodic theorem (see \cite{Oseledet_MSR}) states that there exists numbers $\lambda_1 \leq\ldots\leq \lambda_d$ such that for $\mu$-almost every point $x$ and every vector $v$ in the tangent space at $x$, the limit
\begin{equation}
\underset{N\to\infty}{\lim}\frac{\ln\|DT^N(x)v\|}{N}
\end{equation}
exists and equals one among $\lambda_1,\ldots,\lambda_d$. These numbers are the \textbf{Lyapunov exponents} of the map. For the rest of the section, we will assume that we have a quasiperiodic orbit $(x_n)$ filling out a 2-dimensional torus $\mathbb{T}^2$ and therefore $d=2$ and $M=\mathbb{T}^2$. 

One of the properties of the largest Lyapunov exponent $\lambda_2$ is that for almost every vector $v$ in the tangent space at $x_0$, the quantity $N^{-1}\ln\|DT^{N}(x_0)v\|$ converges to $\lambda_2$. Note that this limit can be expressed as an average along the trajectory $(x_n)$ in the following manner:
\[\frac{\ln\|DT^n(x_0)v\|}{N}=\sum_{n=0}^{N-1}\frac{\ln\|v_n\|}{N}, \mbox{ where }v_n=DT(x_{n-1})\frac{v_{n-1}}{\|v_{n-1}\|}, v_0=v. \]
We obtain the same limit by taking our Weighted Birkhoff average instead of a uniformly weighted average as above. We observe (without proof) that for quasiperiodic orbits, we get super-convergence to the same limit by taking a Weighted Birkhoff average in the following manner:
\[\lambda_2=\underset{N\to\infty}{\lim}\sum_{n=0}^{N-1}\hat{w}_{n,N}\ \ln\|v_n\|,  \mbox{ where }v_n=DT(x_{n-1})\frac{v_{n-1}}{\|v_{n-1}\|}, v_0=v.\]
See Chapter 3 from \cite{ChaosBook} for an explanation of this method.

Once we have calculated $\lambda_2$, $\lambda_1$ can be obtained from the the sum $\lambda_1+\lambda_2$ of the Lyapunov exponents.
$\lambda_1+\lambda_2$ can be expressed as the Birkhoff average of the function $f_1(x):=\ln|\det DT(x)|$. 
Note that $\det DT(x_n) > 0$ for an orientation preserving map $T$. 
Therefore, we have 
\[\lambda_1+\lambda_2=\underset{N\to\infty}{\lim}\frac{\ln|\det DT^N(x)|}{N}=\underset{N\to\infty}{\lim}B_N(f_1).\]
As mentioned before, the quantity $\mbox{B}_N(f_1)$ converges to  the limit $\int_X\ln|\det DT(x)|d\mu(x)$ and we can obtain super-convergence to the same limit using the quantity $\Q_N(f_1)$ as written below:
\begin{equation} \label{eqlyapsum}
\lambda_1+\lambda_2=\underset{N\to\infty}{\lim}\Q_N(f_1)=\underset{N\to\infty}{\lim}\sum\limits_{n=0}^{N-1}\hat{w}_{n,N}\ln|\det DT(x_n)|.
\end{equation}

\subsection{Related methods}\label{Relatedmethods}
See \cite{Laskar99,Laskar93a,Laskar93b,Laskar03,seara:villanueva:06,luque:villanueva:14} for references to earlier methods for computing Birkhoff averages along a quasiperiodic orbit. 

For higher dimensional quasiperiodicity ($d>1$) Laskar  \cite{Laskar99} has an interesting technique of finding the frequency $\sigma$ that maximizes a function $\phi(\sigma)$. That corresponds to our rotation rate $\rho$.  Each evaluation of  $\phi(\sigma)$ requires application of the window filter (whereas our method uses only one application of $\Q$). But his has the advantage of subtracting off this frequency and repeating this method to find the next frequency. Hence he has an automatic method for finding multiple rotation rates. Presumably our method could be combined with his to improve on both methods -- when dealing with higher dimensional quasiperiodicity.

A. Luque and J. Villanueva \cite{seara:villanueva:06,luque:villanueva:14}
develop fast methods for obtaining rotation numbers for
smooth or analytic functions on a quasiperiodic torus, sometimes with
quasiperiodic forcing with several rotation numbers. 
The paper \cite{luque:villanueva:14} develops a technique to compute rotation numbers (but not other function integrals) with error satisfying $|error| \le C_p N^{-p}$ for $N\geq 2^p$ where $C_p$ is a constant. 
Their method of computation is defined  recursively, with the $p+1^{st}$ method being defined in terms of the $p^{th}$. 
As $p$ increases the computational complexity increases for fixed $N$. 
If $\tau(p,N_p)$ is their computation time when using $N=N_p$, it appears that $\tau(p,N_p)/N_p\rightarrow\infty$ as $p\rightarrow\infty$ (and $N_p \to \infty)$. 
In comparison, computation time for our Weighted Birkhoff average is simply proportional to $N$ since it requires a sum of $N$ numbers. Fig. 11 from \cite{luque:villanueva:14} shows the rate of convergence to the rotation number for a quasiperiodic orbit arising from the R3BP.  
There, they get 30-digit accuracy for the rotation number using approximately $2,000,000$ trajectory points while we get the same accuracy with $20,000$. The rate of convergence of their method is $\approx N^{-7.8}$. 
Their methods were extended in \cite{luque:RotDeriv} to compute the derivatives of the rotation numbers for 1-parameter families of circle diffeomorphisms, and in \cite{luque:fourier} to compute the Fourier coefficients of the conjugacy function (i.e. the change of variables between 
the map and rigid rotation).

{\bf Newton methods in the literature.}
An alternative approach to our approach for finding a conjugacy is considered in~\cite{Haro, Cabre, Llave05, Huguet}. In the current paper, we  are assuming that we are starting with only a set of iterates for a single finite length forward trajectory, rather than having access to the functional form of the defining equation. In contrast, the approach in the papers above assumes access to the full form of the original defining equations. 
The Fourier series for the conjugacy is obtained and validated by using automatic differentiation. In addition, in this current paper, we assume that we start with a point in an invariant torus with quasiperiodic dynamics, whereas the methods referenced above include a fast Newton's method for finding invariant tori, Lyapunov multipliers, and invariant stable and unstable manifolds. See also Jorba~\cite{NumericQuasi5}.

Several variants of the Newton's method have been employed to determine quasiperiodic trajectories in different settings. In \cite{NumericQuasi3} a variant of Newton's method was applied to locate periodic or quasi-periodic relative satellite motion in a non-linear, non-conservative setting. 
 
A PDE-based approach was taken in \cite{NumericQuasi4}, where the authors defined an invariance equation which involves partial derivatives. The invariant tori are then computed using finite element methods. See also Section 2 in \cite{NumericQuasi4} for more references on the numerical computation of invariant tori.
In \cite{luque:KAM}, the authors used an application of the Newton's method to compute elliptical, low dimensional invariant tori in Hamiltonian systems.

Computing the conjugacy to the rigid rotation is key to the methods of de la Llave et al \cite{Haro, Cabre, Llave05, Huguet}. They use this to obtain a numerical “proof” of the existence of tori, and, what they call “a posteriori KAM theory”.


\section{Why our method works and when $N$ must be large}\label{sect:nearly-rational}
We will assume that we have a $C^\infty$ quasiperiodic map $F:\torus\to \torus$ and a $C^M$ map $f:\torus\to\mathbb{C}$. 
To better understand convergence on an averaging method applied to an $f$, write
\begin{equation}
f(\theta) = \sum_k a_kf_k(\theta), \mbox { where } f_k(\theta) := e^{2\pi i k\cdot\theta}.
\end{equation}
In particular $a_0 = \int_{\torus}f.$
$\Q_N$  is of course a weighted average of $(f(n\rho))_{n=0}^{N-1}$ and by Theorem~\ref{thm:A}  $\Q_N(f)\to a_0$ as $N\to\infty$. 
Since  an averaging process is linear, 
we can define $\psi_{N,k,\rho} = \Q_N(f_k) $, a collection of numbers that is independent of the choice of $f$. (More generally one could define $\psi^{w}_{N,k,\rho}$ for any weighting function $w$ and the advantages of different choices of $w$ are reflected in the magnitudes $|\psi^{w}_{N,k,\rho}|$). They depend only on the averaging method (e.g., $\Q_N$), $k$, $N$, and the rotation number $\rho$. This set encapsulates all errors that arise in the use of $\Q_N$. In particular,
$f(\theta) = \sum_k a_kf_k(\theta)$,
we have
\begin{equation}
\Q_N(f) = 
a_0+ \sum_{k\ne 0} a_k \psi_{N,k,\rho}. 
\end{equation}
In particular $\psi_{N,0,\rho} = 1$ and for each $k\ne 0$, $\lim_{N\to\infty}\psi_{N,k,\rho}=0$.
The rate of convergence to $0$ depends on $\rho$ and $k$ and convergence can be slow when $e^{2\pi i k\cdot\rho}-1 \approx 0$, as we show in the next section.

\begin{figure}[]
\centering

\subfigure[ ]{\includegraphics[width = .45\textwidth,height=.45\textwidth]{\Path 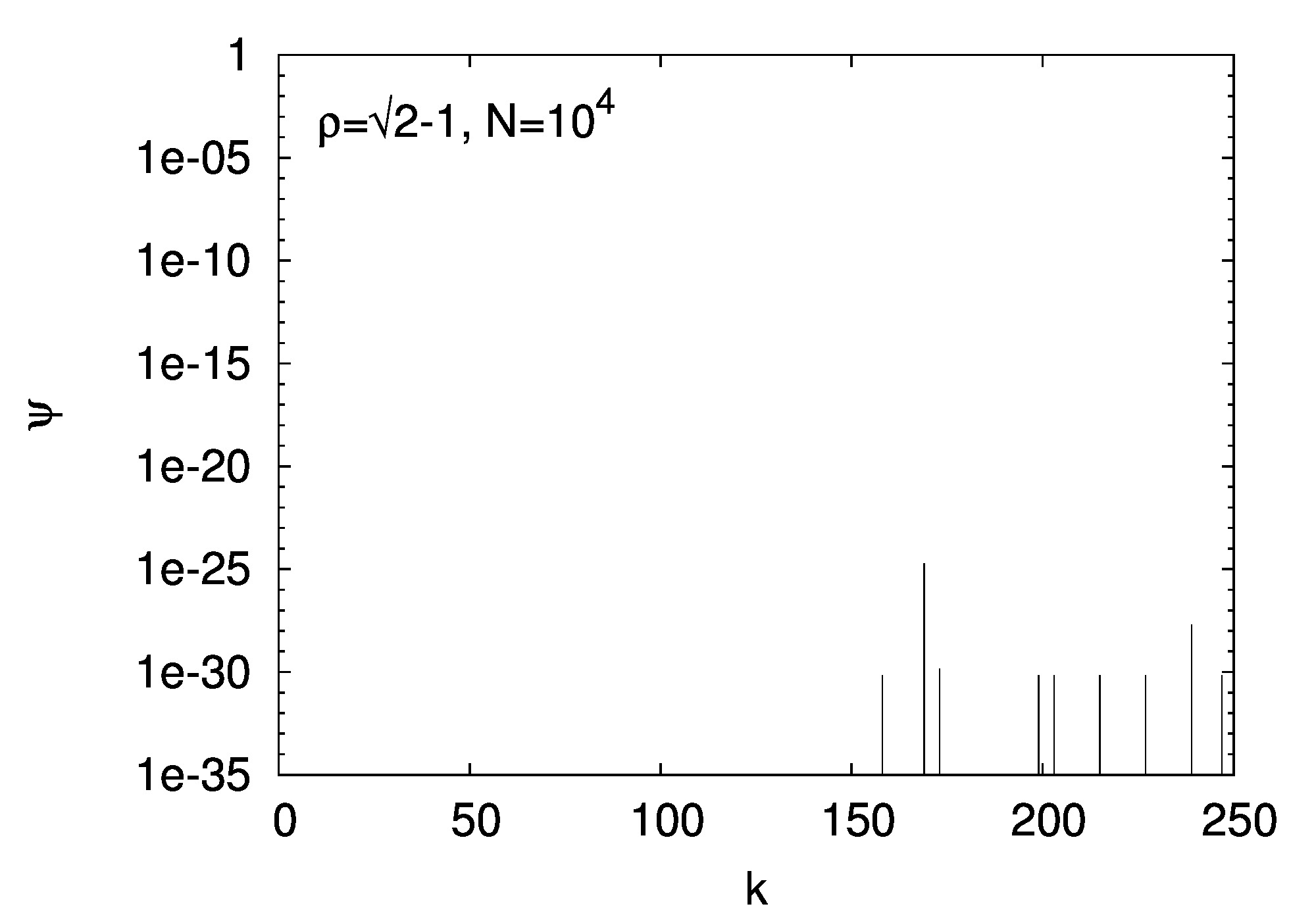}}
\subfigure[ ]{\includegraphics[width = .45\textwidth,height=.45\textwidth]{\Path 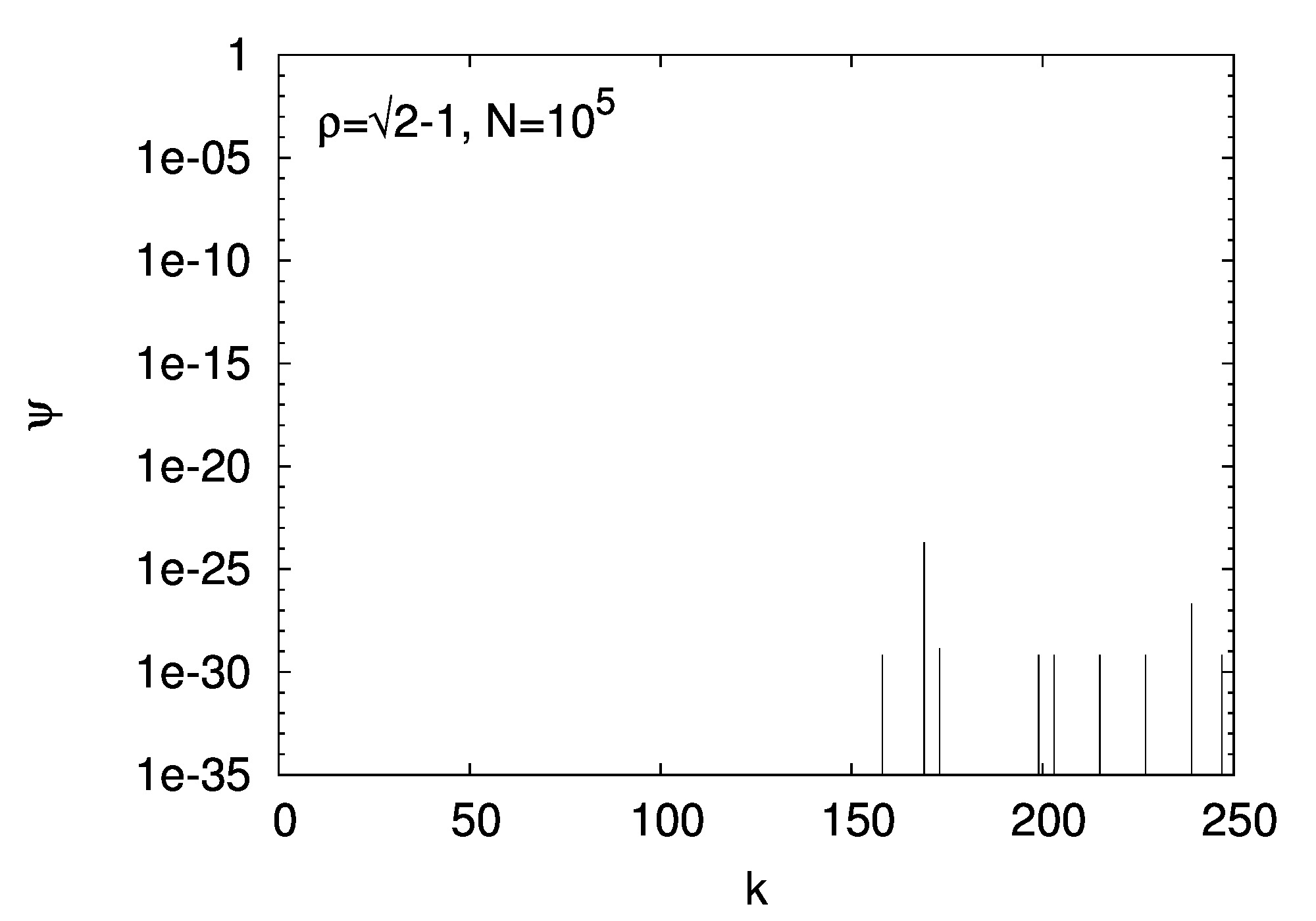}}
\subfigure[ ]{\includegraphics[width = .45\textwidth,height=.45\textwidth]{\Path 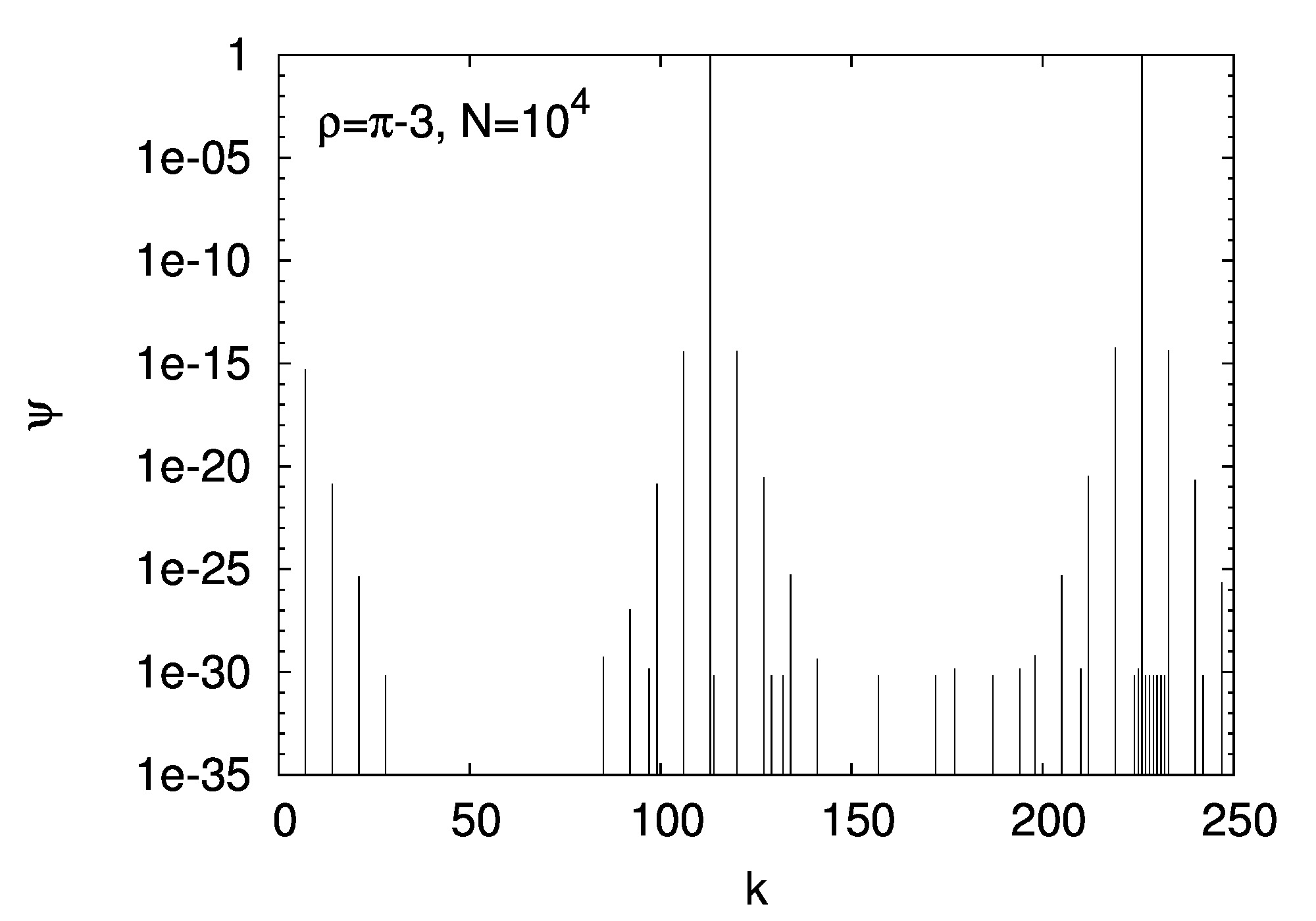}}
\subfigure[ ]{\includegraphics[width = .45\textwidth,height=.45\textwidth]{\Path 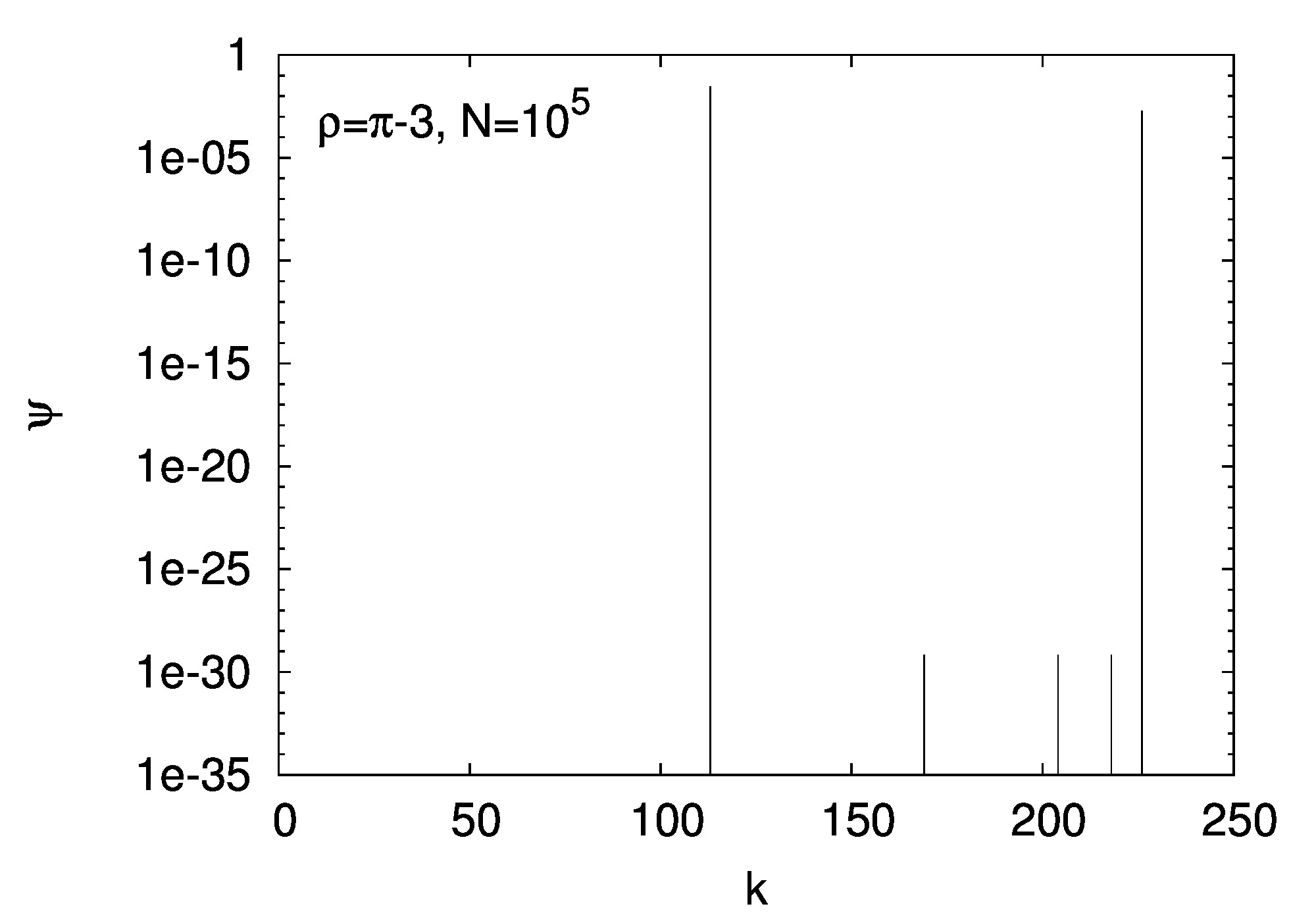}}
\caption{\textbf{Estimations of the numerical error in the calculations of Fourier coefficients.}
Panels (a) and (b) show the value of $\psi_{k,N,\rho}=\Q_{N}(f_k)$ for $\rho=\sqrt{2}-1$, when the number of iterations is $N=10^4$ and $10^5$, respectively. Panels (c) and (d) show for $\rho=\pi-3$, when $N=10^4$ and $10^5$, respectively.  }
\label{fig:error-estimate-psi}
\end{figure}
To investigate the effects of the rotation number $\rho$ on the numerical errors in calculating a Fourier coefficient, 
we show $\psi_{N,k,\rho}=\Q_{N}(f_k)$ for each $k~(0\le k \le 250)$ for two cases of $\rho=\sqrt{2}-1$ and $\pi-3$ in Fig. \ref{fig:error-estimate-psi}. 
The figure suggests that $N=10^5$ is sufficient to calculate Fourier coefficients 
for $\rho=\sqrt{2}-1$, whereas for  $\rho=\pi-3$ even $N=10^5$ is not sufficient.

\subsection{Sketch of proof  of Theorem \ref{thm:A}} \label{sec:proof}
Here we sketch a proof that enables us to determine what happens when a rotation number is near a rational number.

Note that for any constant $a_0$, $\Q_N(a_0 + f)=a_0+ \Q_N(f)$
so for simplicity we will assume $f$ has mean $a_0=0$.
Let \boldmath $I$ \unboldmath denote the identity operator and \boldmath $U$ \unboldmath denote the Koopman operator on $L^2(X,\mu)$, defined as 
\[(Uf):x\mapsto f(F(x)), \mbox{ for all }f\in L^2(X,\mu).\]

The idea of the proof is to provide two different estimates of the quantity
\begin{equation}\label{eqn:Q_N_m}
\Q_N(U-I)^m f = \frac{1}{N}\sum_{n} w\left(\frac{n}{N}\right) (U-I)^m f(n\rho).
\end{equation}

First let $m=1$ in Eq. \ref{eqn:Q_N_m}. Then
\[\begin{split}
\frac{1}{N}\sum_{n} w\left(\frac{n}{N}\right) (U-I) f(n\rho) &= \frac{1}{N}\sum_{n} w\left(\frac{n}{N}\right) f((n+1)\rho) - \frac{1}{N}\sum_{n} w\left(\frac{n}{N}\right) f(n\rho) \\
&= \frac{1}{N}\sum_{n} \left[w\left(\frac{n+1}{N}\right) - w\left(\frac{n}{N}\right)\right] f(n\rho) \\
\end{split}\]
Now taking absolute values on both sides give,
\[\left| \frac{1}{N}\sum_{n} w\left(\frac{n}{N}\right) (U-I) f(n\rho) \right| \leq \frac{1}{N}\|f\|_{C^0}\sum_{n} \left| \frac{1}{N}w^{(1)}\left(\frac{n}{N}\right) \right| \leq \frac{1}{N}\|f\|_{C^0} \|w^{(1)}\|_{C^0}.\]
Applying this procedure $m$ times gives a constant $c_m>0$ such that,
\begin{equation}\label{eqn:QN_U}
\left| \Q_N(U-I)^m f \right| \leq c_m N^{-m}\|f\|_{C^0},
\end{equation}
where $c_m$ depends only on the first $m$ derivatives of $w$.

The second way to evaluate Eq. \ref{eqn:Q_N_m} is by using the Fourier series for $f$.
Note that $f$ being 
$C^M$ implies that 
there is a constant $C_{f,m}>0$ such that $|a_k|\leq C_{f,m} \|k\|^{-M}$. If $k\neq 0$, then
\[\Q_N(U-I)^m f_k  = (e^{2\pi i k\cdot\rho}-1)^m \Q_N f_k. \]
Hence, by Eq. \ref{eqn:QN_U}, 
\begin{equation} \label{ineq:psi}
 |\psi_{N,k}|=\left|\Q_N f_k \right| \leq c_m \left(\frac{1}{N \left| e^{2\pi i k\cdot\rho}-1 \right|} \right)^m.
 \end{equation}
This implies
\begin{align}\label{eqn:dio}
\left|\Q_N (f) \right| &= \left| \sum_{k\neq 0} a_k \Q_N(f_k) \right| 
\leq
\sum_{k\neq 0} |a_k| \left| \Q_N(f_k) \right| \notag
\\
& \leq c_m \left(\frac{1}{N}\right)^{m} \sum_{k\neq 0} |a_k|
\left| \frac{1}{e^{2\pi i k\cdot\rho}-1 }\right|^m. 
\end{align}
The claim of the Theorem now follows from Eq. \ref{eqn:Dioph} and the decay rate of $|a_k|$.

\subsection{Difficulties when $\rho$ is approximately rational.}\label{sec:almost_rational}
While Theorem \ref{thm:A} requires $\rho$ to be Diophantine to get fast convergence of $\Q_N(f)$ to $\int f d\mu$ (the integral in the Birkhoff Ergodic Theorem), our computations have finite precision so we have to ask what this condition means in a finite precision world. 
With this in mind, we can restate the Diophantine condition, Ineq. \ref{eqn:Dioph}, saying $\rho$ is {\bf Diophantine} (of class $\beta$) if there is some $\beta>0$ for which 
$$ 
\liminf_{|k|\to\infty}  |k|^{d+\beta} \left| e^{2\pi i k\cdot\rho}-1 \right|>0.
$$
In other words, the values of which $k$ for which $|k|^{d} \left| e^{2\pi i k\cdot\rho}-1 \right|$ should not go to $0$ too fast. This idea suggests defining 
\[
\Delta(k,\rho):= |k|^d \left| e^{2\pi i k\cdot\rho}-1 \right|.
\]For fast convergence in computing Fourier series coefficients $a_k$, the quantity $\Delta(k,\rho)$ should not be too small for the relevant $k$, those for which  $|a_k|$ is likely to be larger than our error threshold, which in this paper is about $10^{-30}$. 

\begin{figure}
\centering
\subfigure[ ]{\includegraphics[width = .46\textwidth] {\Path 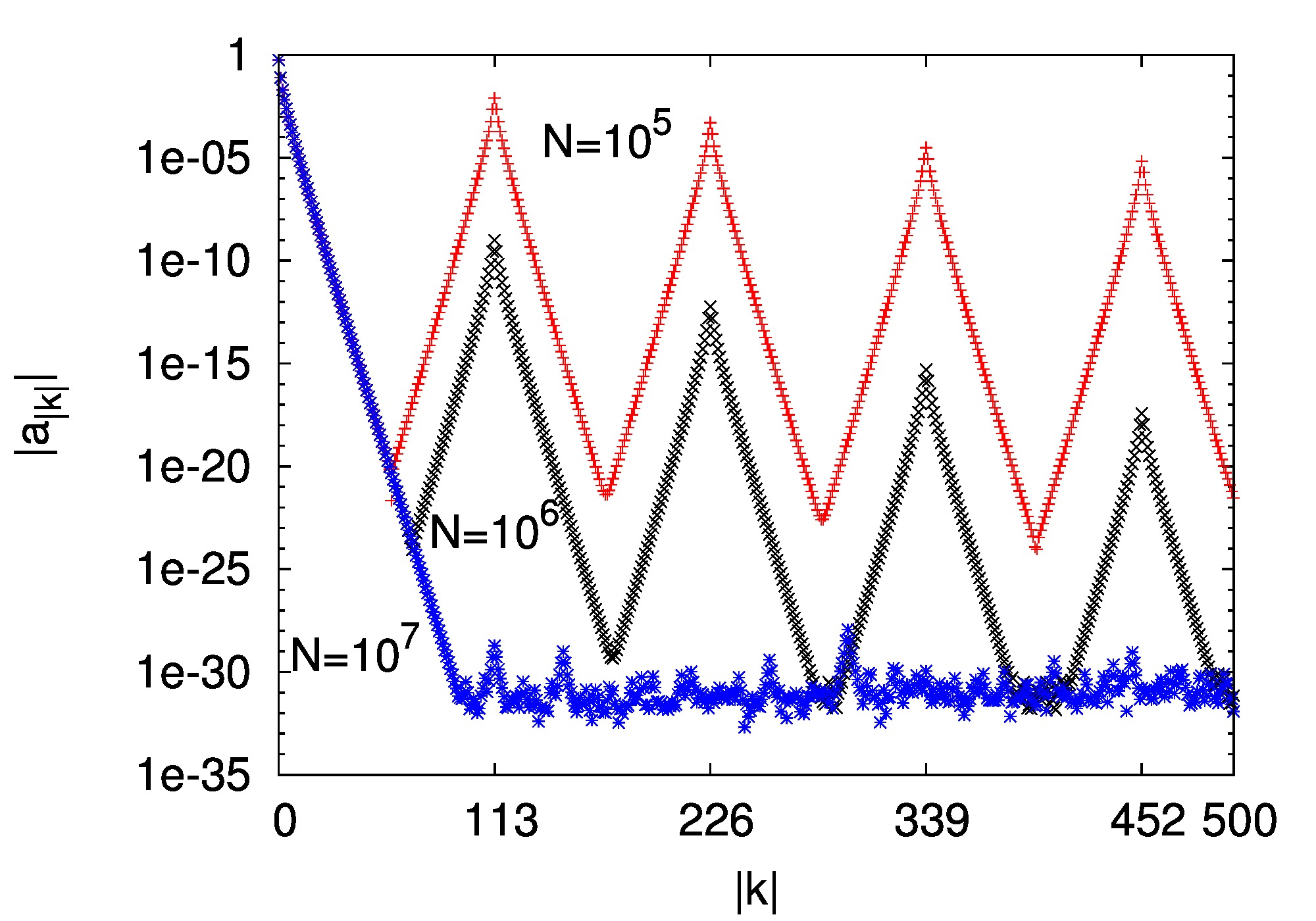}}
\subfigure[ ]{\includegraphics[width = .46\textwidth] {\Path 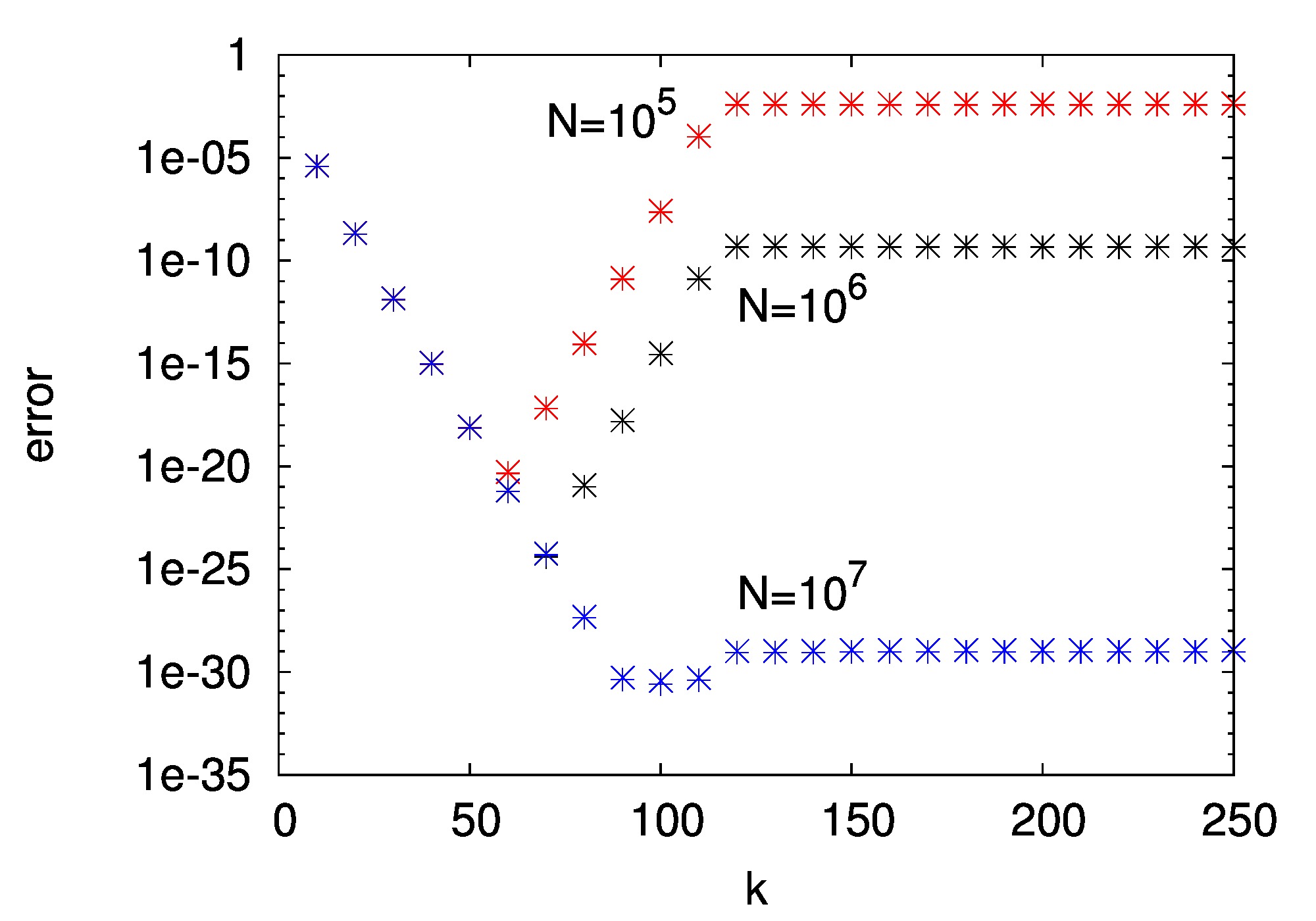}}
\caption{\textbf{When $\rho$ is near a rational.} This figure is for the same case as  Fig. \ref{fig:circle_rotation}(a), but with the rotation number $\rho=\pi-3$. This rotation number is chosen because we require number of iterates $N$ to be large to get accurate results. 
These plots show a slower  convergence compared to the case in $\rho=\sqrt{2}-1$. 
The number of iterates $N$ used to calculate Fourier coefficients is changed from the top to the bottom ($N=10^5$ (top), $10^6$ (middle), $10^7$ (bottom)). Panel (a) shows the norm of the Fourier coefficients $a_k$ of the periodic part of the conjugacy map -- as a function of index $k$. 
There are peaks at the multiples of 113 for $N=10^5$ and $10^6$.
Panel (b) shows the agreement and the disagreement of the Fourier series and the original function in $L^1$ ($+$) and $L^2$ ($\times$) norms in Eq. \ref{deltas} for the corresponding $N$. Note that the $L^1$ and $L^2$ errors are so similar that the 
$+$ and $\times$ overlap, yielding 8-legged spiders. The graph suggests that $N=10^7$ is large enough to calculate Fourier coefficients, whereas for $N=10^5$ and $10^6$, it is best to include only 60 or 70 coefficients respectively, stopping when the error is minimum. }
\label{badrotation}
\end{figure}

{\bf What values of $\Delta(k,\rho)$ are we likely to encounter?}
For $\rho_1 := (\sqrt{5} +1)/2$, the golden mean, we find
\[
\min_{k=2,\cdots,10^6}\Delta(k,\rho_1) \approx 2.655.
\]
It appears that $\displaystyle\liminf_{|k|\to\infty}\Delta(k,\rho_1) = 2.809925\cdots.$
Note the last term in Eq. \ref{eqn:dio} is similar to $\Delta$ except for the leading $|k|$. Of course the power $-m$ can greatly multiply the problem of 
$\left| e^{2\pi i k\cdot\rho}-1 \right|$ being small.
To offset $\left| e^{2\pi i k\cdot\rho}-1 \right|$ being smaller by a factor of 100, we might expect that convergence would require
$N$ in Eq. \ref{eqn:dio} to be larger by a factor of $100$. Of course both are raised to the same power $m$ in these equations. Indeed in Fig. \ref{badrotation} we must increase N by a factor of $100$ to get 30-digit convergence.

An illustration of the problem can be seen for $\rho_2 = \pi-3$
since $\pi \approx \frac{355}{113},$ and $|k|=113$ yields the rather small value  $\Delta(\pm113,\rho_2)\approx 0.021.$ Also   $\Delta(\pm226,\rho_2)\approx 0.085.$
 $\Delta(\pm339,\rho_2)\approx 0.193.$
These suggest slow convergence for those $k$ values. 
When computing a Fourier series for $f$ using $\Q_N$ with this $\rho$,
 when $N < 10^7$, we obtain poor values for many Fourier coefficients as we illustrate in Fig. \ref{badrotation}(a).

{\bf Comparing a function to its Fourier series.}
To estimate how accurate a computed Fourier series of a function $f: \torus \to   \R^1$ is, define  $a_k = \Q_N(f f_{-k})$ 
and  $f_k(\theta) := e^{2\pi i k\cdot\theta}$\ (where $k\cdot\theta$ denotes an inner product in $\R^d$), and 
$$f^{K,N}(\theta) := \sum_{|k| \le K}a_k f_k(\theta).$$
To test how similar $f$ and $f^{K,N}$ are in the $L^2$ and $L^1$ sense, we compute the errors 
\begin{equation}\label{deltas}
\delta^{K,N}_2(f):=\sqrt{\Q_N((f - f^{K,N})^2)} \mbox{ and }\delta^{K,N}_1(f):= \Q_N(|f - f^{K,N}|).
\end{equation}
They are calculated in Fig. \ref{badrotation}(b). 

{\bf A saw-tooth pattern of errors.}
In the above example with $N=10^5$, 
Fig. \ref{badrotation} shows a saw-tooth pattern with peaks at multiples of $113$ where the slopes are all the same. 
Many coefficients have a big error, and all errors here are the result of $\psi_{N,k}$ when $k$ is a non-zero multiple of $113$.
To explain this saw-tooth effect, we note that we have found for this example that $\psi_{N,k}\approx 0$ for $N=10^5$, when
$k$ is not a multiple of $113$ and in particular $k\ne 0$, and here we assume all those values are indeed $0$ to simplify computation except for $k=0$ and $|k| = k^*$ for some $k^* \gg 1$ such as $k^* = 113$.

Suppose we wish to compute the Fourier coefficients of $f$ and determine how accurate the result is. The computed $k^{th}$ coefficient, denoted $\hat a_k$ is 
\begin{align}
\hat a_{|k|}=\hat a_{\pm k}=\Q_N(f(\theta) e^{\pm 2\pi i k\cdot\theta}) =\sum_m a_m\psi_{N,m\pm k} 
\end{align}
which has significant contributions from $m\pm k=0$ and $|m\pm k|= k^*$. Since $k^* \gg 1$, we can ignore $a_{\pm k}$, and we conclude that for small integers $n$,
\[
\hat a_{|k^*\pm n|} \approx a_{|n|}\psi_{N,k^*}.
\]
On the log-linear plot of the graph, the coefficients $a_k$ are almost linear, so it follows from this equation that the erroneous $\hat a_{|k^*\pm n|}$ has the same slope (the absolute value of the derivative), and that each non-zero multiple of $k^*=113$ has the same error pattern. The heights of the peaks are 
$a_{0}\psi_{N,m k^*}$ at $k=m k^*$ where $m=1,2,\cdots$.

{\bf Remark on Fig. \ref{fig:StdMap}.} The panel (b) of this figure shows a jump in the Fourier series terms at coefficient 
$|a_{506}|$ (i.e., $k=253$ in that figure).  
This does not seem to be a numerical artifact.
In fact $\Delta$ does have a local minimum at $506$ but $\Delta(506,\rho)= 1.03$ is not particularly small for such minimum. Furthermore, changing the number of iterates (as mentioned in the caption) does not change the graph. 

\section{Concluding remarks}\label{sec:conclude}

We have developed a straightforward but effective computational tool for quickly computing a large variety of quantities for quasiperiodic orbits. These quantities include rotation vectors, Fourier reconstruction of conjugacy maps, and in some cases Lyapunov exponents. The methods work well in one and higher dimensions. They are effective using both double and quadruple precision, though we have chosen to do most of our calculations in higher precision to show the full possibilities and quick convergence properties of our method. 

The literature on quasiperiodicity is vast and windowing techniques analogous to ours are often used. But our goals in this paper are limited: to introduce the $C^\infty$ Weighted Birkhoff averages $\Q_N\ (=\Q_N^{[1]})$ and $\Q_N^{[2]}$ as numerically useful tools and to present some of its applications. 

We note that the computational time for computing our weighting functions $w(t)$ is almost the same as for the weighting function $w_{sin^2}$ that uses $\sin^2(\pi t)$. Both are equally easy to program. But convergence is far faster with the Weighted Birkhoff averages, as seen in Figs. \ref{fig:mid_circle2},  \ref {fig:2DMap_overview}(b), \ref {fig:circle_rotation}(b), and \ref {fig:2DMap_lyap}.

Quasiperiodic orbits can occur in many different situations, for example, subject to periodic forcing (see Luque and Villanueva \cite{luque:villanueva:14}); as high-dimensional tori that are not simply embedded (see Medvedev et al.~\cite{medvedev_lagrangian_2015}); in the presence of noise; etc. The question of whether our methods extend to these situations  is worthy of further consideration.

{\bf When must $N$ be large to get convergence?}
We have developed some diagnostics in Section \ref{sec:almost_rational} to detect when $N$ must be chosen especially large to get high accuracy -- at least for the $d=1$ dimensional cases.
Computation of $\Delta$ can be used to detect cases when $N$ must be large to get accurate values for Fourier coefficients.  For example we found that because $\Delta(\pm113,\rho_2)$ is so small, $N$ must be increased by a factor of $100$ to get an accurate Fourier series for the conjugacy map. 
One might ask if there are other $k$ for which $\Delta(\pm113,\rho_2)$ is quite small. We find
$\Delta(\pm113,\rho_2) < \Delta(k,\rho_2)$ for all $k\ne 113$ and $k < 10^7$.

\bigskip 

{\bf Acknowledgments:} 
We would like to thank Miguel Sanju\'an and the referees for many helpful suggestions.
YS was partially supported by JSPS KAKENHI grant 17K05360 and JST PRESTO grant JMPJPR16E5.
ES was partially supported by NSF grant DMS-1407087. 
JAY was partially supported by National Research Initiative Competitive grants 2009-35205-05209 and 2008-04049 from the U.S.D.A.

\bibliographystyle{unsrt}
\bibliography{Weighted_calc_bibliography,1D_bibliography}
\end{document}